\documentclass[fleqn,10pt]{wlscirep}
\usepackage[utf8]{inputenc}
\usepackage[T1]{fontenc}
\usepackage{float}
\usepackage{multirow}
\usepackage{subcaption}
\usepackage{graphicx}
\usepackage{changepage}
\usepackage{CJKutf8}
\usepackage{tabularx}
\usepackage{lineno}
\title{Time-Delayed Game Strategy Analysis Among Japan, Other Nations, and the International Atomic Energy Agency in the Context of Fukushima Nuclear Wastewater Discharge Decision}

\author[1]{Mingyang Li}
\author[2]{Han Pengsihua}
\author[1]{Fujiao Meng}
\author[1,*]{Zejun Wang}
\author[1,3]{Weian Liu}
\affil[1]{School of Liberal Arts and Sciences, China University of Petroleum-Beijing at Karamay, Karamay 834000,Xinjiang, China.}
\affil[2]{School of Business Administration, China University of Petroleum-Beijing at Karamay, Karamay 834000,Xinjiang, China.}
\affil[3]{School of Mathematics and Statistics, Wuhan University, Wuhan, Hubei 430072, China.}

\affil[*]{wangzj@cupk.edu.cn (Z.W.)}


\begin{abstract}
On August 24, 2023, the Japanese government made the decision to release treated nuclear wastewater from the Fukushima Daiichi Nuclear Power Plant into the sea. This decision has sparked widespread attention and discussions on an international scale. This paper systematically investigates the time-delayed game strategies between Japan, other countries, and the International Atomic Energy Agency (IAEA) in this context.
This study begins by introducing a payoff matrix and establishing replicator dynamic equations. Taking into account the lag in decision-making in real-world scenarios, time-delay elements are introduced to make the equations more reflective of practical game situations.
Through an analysis of the characteristic roots of the linearized system, this paper delves into the stability of tripartite strategies and the conditions and possibilities of reaching different asymptotically stable states.
Numerical simulations further confirm the practicality and accuracy of this model, particularly in revealing the impact of key parameters on the evolution of the game. The research findings indicate that time delay significantly impacts the stability of decision-making and the evolution trajectories, particularly concerning the strategies of the Japanese government and the International Atomic Energy Agency (IAEA) in the context of nuclear wastewater disposal.
Furthermore, the study underscores the critical importance of scientifically efficient nuclear wastewater treatment technology and highlights the influence of export tax revenue losses on Japan's nuclear wastewater treatment strategies and the significance of international cooperation in addressing this issue.
The innovation of this research lies in its incorporation of time-delay elements from ocean dynamics and government decision-making, providing a solid theoretical foundation for finding optimal strategies for tripartite game entities with delays in the real world.
\end{abstract}
\begin{document}
\begin{CJK}{UTF8}{gbsn}
\flushbottom
\maketitle
%
%
\thispagestyle{empty}

\section{Introduction}

On August 24, 2023, the Japanese government officially decided to release treated radioactive water from the Fukushima Daiichi Nuclear Power Plant into the sea. This decision has sparked widespread international attention and controversy. The discharge of radioactive water from Fukushima is not just a domestic issue for Japan, but a significant public concern affecting marine environments and human health. Previously, Russia set a deadline for Japan to provide data related to the contamination levels of the radioactive water, threatening a ban on the import of Japanese seafood if not complied. In addition, nations such as North Korea and the Solomon Islands have openly criticized Japan's decision to release the contaminated water into the sea. Within Japan, 60\% of the Fukushima residents believe that the water discharge lacks openness and transparency\cite{r3}. Allowing Japan to continue releasing hazardous substances into the ocean could inflict serious damage to the global environment\cite{r4}. Following the commencement of the radioactive water discharge by Japan, various environmental and civic groups in South Korea have spoken out, revealing the dangers of the release and condemning the Japanese government's decision. The outcry from South Korean civic and environmental groups reflects widespread public concern and apprehension. They call on the international community to act collectively in addressing the Fukushima radioactive water issue\cite{r5}. China has firmly called for the establishment of a long-term and effective international monitoring mechanism to ensure the impacts of the radioactive water on the environment and human health are effectively monitored and assessed\cite{r6,r7}. As the intensity of the radioactive water discharge increases, its impact assessment and countermeasures have become matters of broad international concern.

Scholars warn that this could result in long-term and irreversible effects on marine ecosystems, with related industries such as fisheries and the food sector experiencing adverse consequences\cite{r14,1,2}. Firstly, contamination from radioactive materials directly impacts fisheries by rendering fish inedible\cite{r8,r9,r10,r11,r12}. Furthermore, tourism industries associated with the ocean are also expected to suffer\cite{r13}. The release of radioactive water by Japan has multifaceted implications, with coastal countries and regions facing potential threats\cite{3,t1,t2,t3}. On October 5, 2023, the Japanese government initiated its second phase of radioactive water discharge\cite{4}. Notably, prior to this official release, Japan's decision had already evoked widespread responses from the international community, especially neighboring countries like China, which strongly oppose Japan's unilateral decision\cite{5}. The role of the International Atomic Energy Agency (IAEA)—a primary international organization promoting and overseeing the peaceful use of nuclear energy\cite{6}—in the Fukushima radioactive water treatment incident has attracted wide international attention\cite{7}. The global community expects the IAEA to play a proactive role in ensuring that the decision-making process is fair, transparent, and serves both human and environmental interests.

In this context, this study aims to explore the strategic interplay among Japan, other countries, and the International Atomic Energy Agency (IAEA) against the backdrop of radioactive water discharge into the sea. Utilizing an evolutionary game theory model, this research considers the potential impact of time delays on decision-making in the actual game process and delves into the evolving strategies of the Japanese government, other countries, and the IAEA. The study primarily focuses on the following four core issues:
\begin{itemize}
	\item How can one accurately construct the payoff matrices for Japan, other countries, and the IAEA, to reflect the genuine interests of each party?
	\item When incorporating delay factors, how can the strategic stability among the three parties be analyzed, and under what conditions is stability likely achieved?
	\item How do varying delays affect the strategic evolutionary trajectories among the three parties, and do these delays alter the conditions for evolutionary stability?
	\item In the context of a complex international backdrop, how can a balance be struck among the interests and concerns of the parties to seek a solution that gains broad acceptance?
\end{itemize}

Incorporating the concept of time delay into evolutionary game models acknowledges the inherent delays in real-world decision-making. This modification specifically addresses the time required by various actors to make decisions, factoring in delays caused by policy development, stakeholder consultations, and strategic planning. These delays significantly affect the dynamics of strategic interactions, making the model more realistic. By integrating these time considerations, the model offers a more accurate portrayal of interactions over time, allowing for a more precise assessment of the stability of equilibrium points in complex strategic scenarios. This approach thus enhances the realism of evolutionary game models by mirroring the tempora aspects of decision-making processes in real-life situations.

This paper begins with Section II, where a systematic review of relevant literature is presented, providing a robust foundation for the theoretical framework and methodology of the entire research.
In Section III, various symbols and definitions used in the study are elaborated in detail, establishing a clear conceptual foundation for subsequent model construction and analysis.
Section IV builds the payoff matrices describing the interests of Japan, other countries, and the IAEA. Based on this, the element of delay is introduced, taking into account the response time of the parties before making decisions.
In Section V, using methods from the mathematical study of delay dynamical systems, the system is linearized, and its characteristic equations and the signs of its eigenvalues are discussed. This analysis uncovers the stability properties of the strategies of the three parties and offers essential insights into how an asymptotically stable state can be achieved.
Section VI employs numerical simulations to empirically test the theoretical model. The reliability and accuracy of the theoretical analysis are verified, and the evolutionary dynamics and trends of the strategies of the three parties under different parameter combinations are further explored.

Overall, this research introduces delay factors into practical applications of evolutionary game theory for the first time. The stability of equilibrium points in delay differential equations is proven, emphasizing its core and vital role in the model. The introduction of time delays not only offers more precise and scientific theoretical support for global marine environment conservation but also provides a new perspective and depth to the strategic analysis among Japan, other countries, and the International Atomic Energy Agency (IAEA).
\section{Literature review}
\subsection{Harmful effects of radioactive water discharge}
On March 11, 2011, a major earthquake off the northeast coast of Japan and the subsequent tsunami led to an accident at the Fukushima Daiichi Nuclear Power Plant (FDNPP). This resulted in the release of radioactive isotopes into the atmosphere, including \textsuperscript{131}I and \textsuperscript{137}Cs\cite{hir}. As reported by the Tokyo Electric Power Holdings (TEPCO) and the Nuclear Regulation Authority (NRA) of Japan, the radioactive substances mainly consist of isotopes like \textsuperscript{3}H, \textsuperscript{14}C, \textsuperscript{134}Cs, and others\cite{tep,nra}. These isotopes are easily absorbed by marine organisms and seafloor sediments (Buesseler, 2020)\cite{bue}.
Studies suggest that if this radioactive water were to be discharged into the Pacific Ocean, its potential risks could persist for hundreds to thousands of years (Albouy et al., 2019)\cite{alb}. Radioactive materials in the wastewater pose significant threats to both the environment and human health (Clifford and Zhang, 1994; Dufresne et al., 2018)\cite{cli,duf}. For instance, \textsuperscript{137}Cs is a highly soluble radioactive isotope in seawater. It carries long-term radiation risks in the environment and, at high concentrations, is easily absorbed by marine organisms (Buesseler, 2020)\cite{bue}. Additionally, the concentration of \textsuperscript{14}C in treated radioactive water remains high (The Korea Times, 2021b)\cite{the}. It has the potential to enter the biosphere and accumulate in marine ecosystems (Williams et al., 2010)\cite{wil}. Co can penetrate the human body, causing cellular damage (Khajeh et al., 2017)\cite{khaj}. Other radioactive substances, like \textsuperscript{90}Sr, can mimic calcium within the human body, increasing risks of osteosarcoma and leukemia, thereby posing a threat to human health (Khani et al., 2012)\cite{khan}.

In summary, the detrimental impacts of radioactive wastewater on the marine environment are long-term, and its negative consequences are hard to quantify. Radioactive isotopes with long half-lives can remain in the ocean for prolonged periods (Men and Deng et al., 2017)\cite{men}, posing significant threats to the natural environment, marine life, and human health.

\subsection{Evolutionary game theory}

Game theory, initially designed for the analysis of economic behaviors, has now expanded not only to the problems of population evolution but also to multiple fields such as economics, management, and environmental science, assisting in addressing issues related to decision-making, strategy, cooperation, and competition. Evolutionary game theory, inspired by biological evolution theories, deviates from the conventional static game theory. Instead of merely focusing on strategy selections, it emphasizes how strategies evolve through repeated games and trial and error (Weibull, 1997)\cite{wei}. Evolutionary game theory provides a theoretical framework for studying the interaction between individual behavior and adaptability. Within this framework, the evolutionary process is perceived as a dynamic strategy selection process, where the fitness of an individual depends not only on its strategy selection but also gets influenced by its interactions with other individuals (Smith, 1973)\cite{smi}.
Contrary to the fully rational participants assumed in traditional game theory, participants in evolutionary game theory are considered to possess bounded rationality. They pursue the maximization of their utility through continuous learning and strategy adjustments (Von Neumann and Morgenstern, 2007)\cite{von}. This process of constant adjustments and learning results in what's termed as a dynamic equilibrium, rather than the static equilibrium in the traditional sense (Gallardo and Marui, 2016)\cite{gal}. The core objective of evolutionary game theory is to identify the best strategy or sequence of decisions in a given conflict to achieve maximum benefits. The competition between strategies is dynamic, where new strategies emerge over time, challenging the existing ones and might replace the dominant strategies at certain periods. This dynamism in strategies led to the concept of Evolutionarily Stable Strategies (ESS), perceived as the optimal strategy that could prevail in strategy competitions (Smith, 1968)\cite{smi1968}.

\subsection{Application of evolutionary game theory in this study}

Evolutionary game theory offers a unique perspective for understanding and analyzing the strategic interactions and choices among the three parties. In recent years, this theory has been widely applied to various environmental and socio-economic issues.
Liu et al. (2021a)\cite{liua} developed a tripartite evolutionary game model specifically to investigate whether Japan chooses to discharge nuclear wastewater, taking into account other neighboring regions and processing costs. Extending this research, Liu et al. (2021b)\cite{liub} employed various game methodologies, including static games, rank-dependent expected utility (RDEU) games, and sequential games, to deeply analyze the equilibrium strategies of the discharging country and the stakeholder countries under different scenarios. Their findings revealed how to balance internal and external interests and make optimal decisions based on emotional states in various game structures.
Ye (2022)\cite{junye} explored the strategies between flag states and port states in the control of ship ballast water discharge. The study analyzed the role of the International Environmental Protection Organization and its impact on the control strategies of flag states and port states. It was found that relying solely on flag state control is ineffective in handling the discharge of radioactive water, while port state control is an effective supplement that helps reduce the risk of marine pollution. Additionally, reasonable assistance from the International Environmental Protection Organization can motivate flag states to adopt more effective control measures.
Zheng (2022)\cite{zhe} established an evolutionary game model that included fishermen, consumers, and the government to explore how government subsidies affect the system's stability. Xin et al. (2022)\cite{xin} focused on how emotions influence the strategic choices of the Japanese government and fishermen, shedding light on the role emotions play in strategy selection and interactions.
Xu et al. (2022)\cite{xul} concentrated on the interactions among the International Atomic Energy Agency, the discharging country, and the Japan Fisheries Association. They observed that the behaviors of these three parties are highly interdependent, especially under the influence of the environmentally-driven International Atomic Energy Agency.

Most existing research concentrates on the ecological threats of the nuclear wastewater discharge by the Japanese government and delves into the relations between the Japanese government and international organizations from an international relations perspective. Yet, no study has deeply explored the interest games among the Japanese government, other countries, and the International Atomic Energy Organization. Crucially, since Japan officially began discharging nuclear wastewater into the Pacific on August 24, 2023, the ocean currents have caused time lags in the arrival of the nuclear wastewater in the marine territories of different countries or regions. Every game player, based on the present impact, requires time to gather more information to formulate and implement policies, thereby causing a time lag effect. By integrating this time lag effect into the evolutionary game equations, we aim to help various stakeholders more precisely predict the long-term impacts of the nuclear wastewater discharge, offering more scientific, rational, and effective strategic suggestions. This will support global efforts to better address the challenges of nuclear wastewater discharge, fostering sustainable marine environmental development.

\section{Methodology}
\subsection{Symbol definitions}

All symbols and their corresponding meanings are detailed and illustrated in Table \ref{not}.

\begin{table}[H]
	\centering
	\begin{tabular}{@{}ll@{}}
		\toprule
		Parameter & Description  \\ \midrule
		\( J \) & Japan \\
		\( C \) & Other countries \\
		\( IAEA \) & International Atomic Energy Agency \\
		\(x\) & Probability that Japan chooses the discharge strategy \\
		\( 1-x \) & Probability that Japan chooses to halt the discharge \\
		\( y \) & Probability of other countries imposing sanctions on Japan's discharge \\
		\(1-y\) & Probability of other countries not imposing sanctions on Japan's discharge \\
		\( z \) & Probability that IAEA opposes Japan's nuclear wastewater discharge \\
		\( 1-z \) & Probability that IAEA supports Japan's nuclear wastewater discharge \\
		\( C_{DJ }\) & Cost for Japan of discharging to the sea \\
		\( C_{SJ} \) & Cost for Japan of storing nuclear wastewater \\
		\( I_J \) & Japan's international reputation \\
		\( T_{RJ} \) & Reduction in export tax revenue for Japan due to discharge \\
		\( C_{MJ} \) & Ocean monitoring cost for Japan under the discharge scenario \\
		\( C_{HJ} \) & Aid received by Japan from other countries in the non-discharge scenario \\
		\( C_{LC} \) & Litigation compensation received by other countries under the discharge scenario \\
		\( C_{SC} \) & Extra cost for other countries in substituting Japanese products under the discharge scenario \\
		\( C_{MC} \) & Ocean monitoring cost for other countries under the discharge scenario \\
		\( C_{II} \) & IAEA's international reputation \\
		\( E_{DI} \) & IAEA's ecological harm metric value due to nuclear wastewater discharge \\
		\( H_{RI} \) & IAEA's health risk metric value due to food contamination \\
		\( C_{MI} \) & Ocean monitoring cost for IAEA under the discharge scenario \\
		
		\bottomrule
	\end{tabular}
	\caption{Description of parameters.}
	\label{not}
\end{table}

\subsection{Basic assumptions}

\begin{enumerate}[label={\textbf{Assumption} \arabic*:},left=0cm]
	\item In this three-party game, Japan (J), other countries (C), and the International Atomic Energy Agency (IAEA) are set as agents. Due to limited information and rationality, each party's decisions are constrained.
	\item Using the evolutionary game theory approach, the above three parties are treated as stakeholders. The dynamic game between these three stakeholders encompasses each one's pursuit of interests and choice of actions, all based on the principle of maximizing benefits.
	\item In this game, Japan (J), being the discharger of nuclear-contaminated water, may weigh its domestic economic benefits against its international reputation and decide whether to continue the discharge. The probability of discharging is denoted by \(x\), and not discharging is \(1-x\).
	\item Countries other than Japan (C) act as other stakeholders. Their attitudes can be classified into sanctioning and not sanctioning. Some countries might express discontent with Japan's sea discharge strategy to safeguard international environmental and marine resources interests, adopting sanctions with a probability of \(y\). Meanwhile, others might opt not to sanction with a probability of \(1-y\).
	\item The International Atomic Energy Agency (IAEA) will consider the environmental and health risks of nuclear wastewater discharge and its own international image. It might oppose Japan's discharge strategy with a probability of \(z\) or support it with a probability of \(1-z\).
	\item This study solely analyzes the game strategies of each party from an economic perspective, without considering the influence of political and other factors on each participant's strategic choice.
	\item To facilitate the research on uniformly distributed scenarios, this study assumes that the time delays experienced by all three parties are equal, denoted by \(\tau\). This assumption is made to eliminate the time delay differences among the three parties, providing a fair benchmark during the gaming process.
	\item The international aid Japan receives in the non-discharge scenario discussed in this paper primarily refers to technological assistance. Such technical assistance mainly includes offering advanced nuclear wastewater treatment technologies, sharing best practices in nuclear wastewater treatment, providing expert teams for on-site guidance and training, and engaging in joint research and development with Japan, among others.
\end{enumerate}

\begin{figure}[H]
	\centering
	\includegraphics[width=0.8\linewidth]{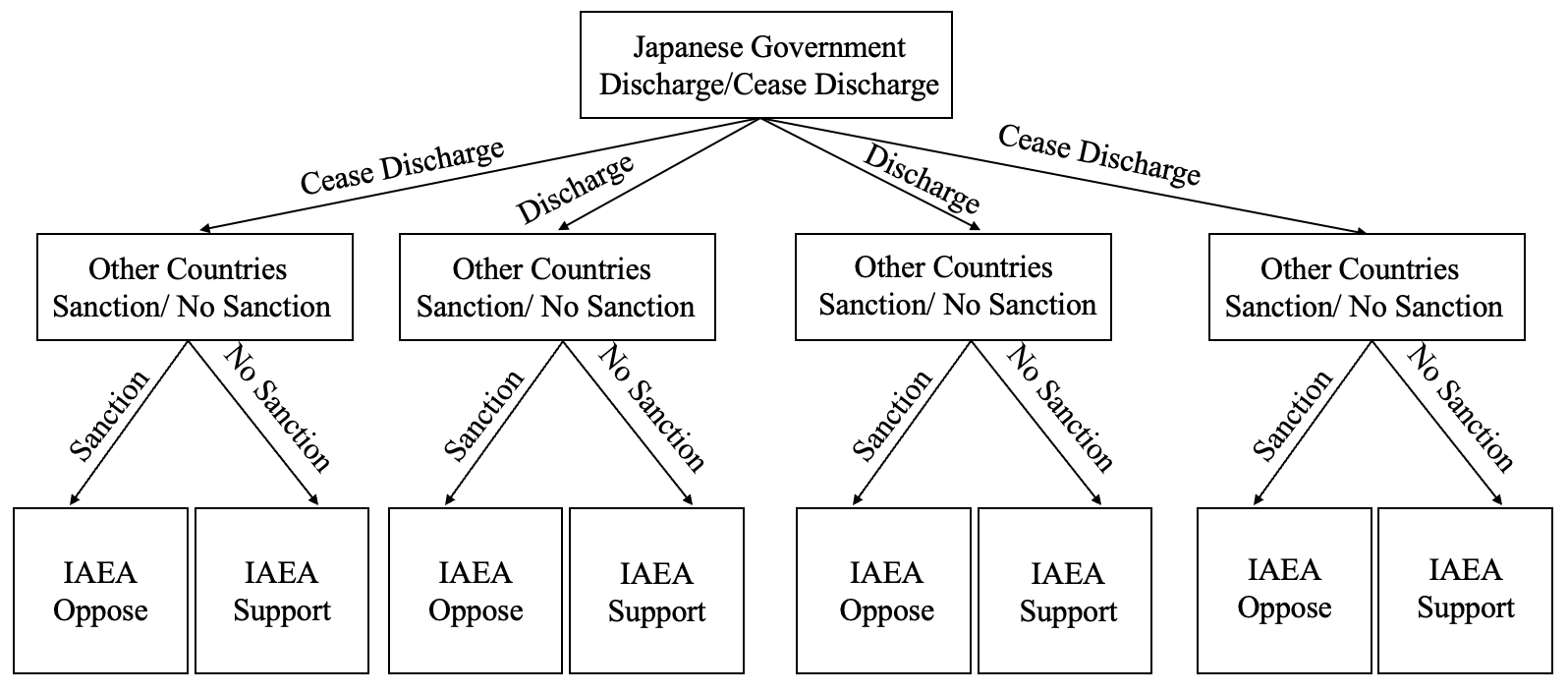}
	\caption{All strategy combinations during the gaming process.}
	\label{fig:screenshot010}
\end{figure}

During the three-party gaming process, all strategy combinations in the game are as shown in Figure~\ref{fig:screenshot010}.

\subsection{Evolutionary game model}

Evolutionary game theory is a mathematical theory that studies how strategies evolve and change over time in a dynamic environment. In international relations, countries and international organizations make decisions based on their payoff matrix. The dynamic nature and mutual influence of these decisions make evolutionary game theory an ideal tool to study these complex interactions.

In the Fukushima nuclear wastewater discharge incident, we particularly focus on the strategic interactions among the Japanese government, other countries, and the International Atomic Energy Agency (IAEA). The decisions of these three parties are not only based on their respective internal interests but are also influenced by the strategic choices of the others. To this end, we have constructed a payoff matrix to describe the game relationship among the three parties.

The payoff matrix includes the Japanese government, other countries, and the IAEA, as shown in Table~\ref{pm}. This matrix illustrates the payoff for one party when the other two choose specific strategies.

	\begin{table}[H]
	\centering
	\small
	\caption{Payoff matrix for Japan, other countries, and IAEA.}
	\resizebox{\textwidth}{!}{ 
		\begin{tabular}{cccccc}  
			\bottomrule
			\multicolumn{4}{c}{\multirow{2}{*}{\textbf{Stakeholders}}} & \multicolumn{2}{c}{\textbf{IAEA}}  \\
			\cline{5-6}
			\multicolumn{4}{c}{} & Oppose $(z)$ & Agree  $(1 - z)$ \\
			\hline
			& \multirow{4}{*}{Discharge $(x)$} & & \multirow{2}{*}{Sanction $(y)$} & $
			\begin{aligned}
				&-I_{J}-C_{L C}-T_{R J}-C_{D J}-C_{M J}, \\
				&-C_{S C}+C_{L C}-C_{M C}, \\
				&-C_{M I}+C_{I I}-E_{D I}-H_{RI}
			\end{aligned}$ & $
			\begin{aligned}
				&-I_{J}-C_{L C}-T_{R J}-C_{D J}-C_{M J}, \\
				&-C_{S C}+C_{L C}-C_{M C}, \\
				&-C_{M I}-E_{D I}-H_{RI}
			\end{aligned}$ \\  
			& & & & &  \\  
			& & & \multirow{2}{*}{No Sanction $(1-y)$} & $
			\begin{aligned}
				&-C_{D J}-C_{M J}, \\
				&-C_{M C}, \\
				&-C_{M I}+C_{I I}-E_{D I}-H_{RI}
			\end{aligned}$ & $
			\begin{aligned}
				&-C_{D J}-C_{M J}, \\
				&-C_{M C}, \\
				&-C_{M I}-E_{D I}-H_{RI}
			\end{aligned}$ \\  
			& & & & &  \\  
			& \multirow{4}{*}{No Discharge $(1-x)$ } & & \multirow{2}{*}{ Sanction$(y)$} & $
			\begin{aligned}
				&C_{HJ}-C_{SJ}, \\
				&-C_{HJ}, \\
				&C_{I I}
			\end{aligned}$ & $
			\begin{aligned}
				&C_{HJ}-C_{S J}, \\
				&-C_{HJ}, \\
				&0
			\end{aligned}$\\  
			& & & & & \\ 
			\multirow{-24.5}{*}{\textbf{Japan}} & & \multirow{-24.5}{*}{\textbf{Other Countries}} & \multirow{2}{*}{No Sanction  $(1-y)$} & $
			\begin{aligned}
				&-C_{S J}, \\
				&0, \\
				&C_{I I}
			\end{aligned}$ &  $
			\begin{aligned}
				&-C_{S J}, \\
				&0, \\
				&0
			\end{aligned}$\\  
			& & & &  &  \\  
			\bottomrule
		\end{tabular}	\label{pm}
	} 
\end{table}
Assuming that the expected utility of the Japanese government discharging nuclear wastewater is represented by \(U_{11}\), the expected utility of not discharging is \(U_{12}\), and the average expected utility is \(\bar{U}_{1}\). Thus, we have:

	\begin{equation}
		\begin{aligned}
			&U_{11}=yz(-I_{J}-C_{LC}-T_{RJ}-C_{DJ}-C_{MJ})+y(1-z)(-I_{J}-C_{LC}-T_{RJ}-C_{DJ}-C_{MJ})\\
			&+(1-y)z(-C_{DJ}-C_{MJ})+(1-y)(1-z)(-C_{DJ}-C_{MJ})\\
			&=-y(I_J+C_{LC}+T_{RJ})-(C_{DJ}+C_{MJ})
		\end{aligned}
	\end{equation}
	
	\begin{equation}
		\begin{aligned}
			&U_{12}=yz(C_{HJ}-C_{SJ})+y(1-z)(C_{HJ}-C_{SJ})\\
			&+(1-y)z(-C_{SJ})+(1-y)(1-z)(-C_{SJ})\\
			&=yC_{HJ}-C_{SJ}
		\end{aligned}
	\end{equation}
	
	\begin{equation}
		\bar{U}_{1}=xU_{11}+(1-x)U_{12}
	\end{equation}
	The replicator dynamics equation for the Japanese government is represented as \(S(x)\), as shown in equation~(\ref{sx}).
	
	\begin{equation}
		\begin{aligned}
			S(x)=\frac{dx}{dt} & =x(U_{11}-\bar{U}_{1})\\
			&=x(1-x)(U_{11}-U_{12}) \\
			& =x(1-x)[-y(I_J+C_{LC}+T_{RJ}+C_{HJ})-(C_{DJ}+C_{MJ}-C_{SJ})]
		\end{aligned}
		\label{sx}
	\end{equation}
	
	Assuming the expected utility of other countries sanctioning the Japanese government is represented by \(U_{21}\), the expected utility of other countries not sanctioning the Japanese government is \(U_{22}\), and the average expected utility is \(\bar{U}_{2}\). Thus, we have:
	
	\begin{equation}
		\begin{aligned}
			&U_{21}=xz(-C_{SC}+C_{LC}-C_{MC})+x(1-z)(-C_{SC}+C_{LC}-C_{MC})\\
			&+(1-x)z(-C_{HJ})+(1-x)(1-z)(-C_{HJ})\\
			&=x(-C_{SC}+C_{LC}-C_{MC}+C_{HJ})-C_{HJ}
		\end{aligned}
	\end{equation}
	
	\begin{equation}
		\begin{aligned}
			&U_{22}=xz(-C_{MC})+x(1-z)(-C_{MC})+(1-x)z\cdot0+(1-x)(1-z)\cdot 0\\
			&=-xC_{MC}
		\end{aligned}
	\end{equation}
	
	\begin{equation}
		\bar{U}_{2}=yU_{21}+(1-y)U_{22}
	\end{equation}
	
	The replicator dynamics equation for other countries is represented as \(G(y)\), as shown in equation~(\ref{gy}).
	
	\begin{equation}
		\begin{aligned}
			G(y)=\frac{dy}{dt} & =y(U_{21}-\bar{U}_{2})=y(1-y)(U_{21}-U_{22}) \\
			& =y(1-y)[x(-C_{SC}+C_{LC}+C_{HJ})-C_{HJ}]
		\end{aligned}
	\label{gy}
	\end{equation}

	Assuming the anticipated utility of IAEA's stance against emissions as \(U_{31}\), while the utility for supporting emissions is represented as \(U_{32}\), and the average expected utility is \(\bar{U}_{3}\). Thus, we have:
		\begin{equation}
			\begin{aligned}
				&U_{31}=xy(-C_{M I}+C_{I I}-E_{D I}-H_{RI})+x(1-y)(-C_{M I}+C_{I I}-E_{D I}-H_{RI})+(1-x)y\cdot C_{I I}+(1-x)(1-y)\cdot C_{I I}\\
				&=-x(C_{MI}+E_{DI}+H_{RI})+C_{II}
			\end{aligned}
		\end{equation}
		
		\begin{equation}
			\begin{aligned}
				&U_{32}=xy(-C_{M I}-E_{D I}-H_{RI})+x(1-y)(-C_{M I}-E_{D I}-H_{RI})+(1-x)y\cdot 0+(1-x)(1-y)\cdot 0\\
				&=-x(C_{M I}+E_{D I}+H_{RI})
			\end{aligned}
		\end{equation}
		
		\begin{equation}
			\bar{U}_{3}=zU_{31}+(1-z)U_{32}
		\end{equation}
	
		The replicator dynamic equation for IAEA, represented as \(P(z)\), is given by equation~(\ref{pz}).
		
		\begin{equation}
			\begin{aligned}
				P(z) & = \frac{dz}{dt} = z(U_{31} - \bar{U}_{3}) = z(1-z)(U_{31} - U_{32}) \\
				& = z(1-z)C_{II}
			\end{aligned}
			\label{pz}
		\end{equation}

		To delve deeper into the stable states within the framework of evolutionary game theory, we integrate equations~(\ref{sx}),~(\ref{gy}), and~(\ref{pz}) to establish the replicator dynamic system depicted in~(\ref{fxgypz}).
		
		\begin{equation}
			\label{fxgypz}
			\left\{
			\begin{aligned}
				S(x) &= \frac{d x}{d t} = x(1-x)[-y(I_J+C_{LC}+T_{RJ}+C_{HJ})-(C_{DJ}+C_{MJ}-C_{SJ})]\\
				G(y) &=\frac{d y}{d t} = y(1-y)[x(-C_{SC}+C_{LC}+C_{HJ})-C_{HJ}]\\
				P(z) &= \frac{d z}{d t} = z(1-z)C_{II}\\
			\end{aligned}
			\right.
		\end{equation}

		In the real world, decisions made by different game entities are not always instantaneous, implying the possibility of a certain time lag. For instance, a nation or organization might require time to assess and react to another country's decisions. Further, factors like ocean currents may cause a delay before nuclear contamination reaches the waters of various nations. To incorporate this time-lag effect, we modify the aforementioned replicator dynamic equations by introducing a fixed time delay \(\tau\). We can describe the system considering this delay as shown in equation~(\ref{delay}).
		
		\begin{equation}
			\label{delay}
			\left\{
			\begin{aligned}
				S(x) &= \frac{d x}{d t}  = x(t-\tau)(1-x(t-\tau))[-y(t-\tau)(I_J+C_{LC}+T_{RJ}+C_{HJ})-(C_{DJ}+C_{MJ}-C_{SJ})]\\
				G(y) &=\frac{d y}{d t} = y(t-\tau)(1-y(t-\tau))[x(t-\tau)(-C_{SC}+C_{LC}+C_{HJ})-C_{HJ}]\\
				P(z) &= \frac{d z}{d t}  = z(t-\tau)(1-z(t-\tau))C_{II}\\
			\end{aligned}
			\right.
		\end{equation}
		
		Where \(\tau\) represents the time-lag parameter, indicating the influence of the state from \(\tau\) time units in the past on the present state.
		
		Based on the replicator dynamic system in equation~(\ref{delay}), the evolutionarily stable states for the Japanese government, other nations, and IAEA can be determined as: \(\gamma_{1}(0,0,0)\), \(\gamma_{2}(1,0,0)\), \(\gamma_{3}(0,1,0)\), \(\gamma_{4}(0,0,1)\), \(\gamma_{5}(1,1,0)\), \(\gamma_{6}(1,0,1)\), \(\gamma_{7}(0,1,1)\), \(\gamma_{8}(1,1,1)\), and \(\gamma_{9}(x^{*}, y^{*}, z^{*})\), where \(\gamma_{9}(x^{*}, y^{*}, z^{*})\) is determined by equation~(\ref{de}).

	\begin{equation}
		\left\{
		\begin{aligned}
			& -y(t-\tau)(I_J+C_{LC}+T_{RJ}+C_{HJ})-(C_{DJ}+C_{MJ}-C_{SJ})=0\\
			& x(t-\tau)(-C_{SC}+C_{LC}+C_{HJ})-C_{HJ}=0\\
			&C_{II}=0\\
		\end{aligned}
		\right.
	\label{de}
	\end{equation}
	
	\section{Evolutionary stability analysis}
	\subsection{Asymptotic stability analysis of the three parties}

	In asymmetric games, the evolutionarily stable strategy is a pure strategy equilibrium, so we only consider the asymptotic stability of pure strategies, hereinafter referred to as stability. Eight equilibrium points have already been obtained in Section 4: \(\gamma_{1}(0,0,0)\), \(\gamma_{2}(1,0,0)\), \(\gamma_{3}(0,1,0)\), \(\gamma_{4}(0,0,1)\), \(\gamma_{5}(1,1,0)\), \(\gamma_{6}(1,0,1)\), \(\gamma_{7}(0,1,1)\), and \(\gamma_{8}(1,1,1)\). The stability of these eight points is analyzed in the appendix, and the results are presented in Table~\ref{st}.

	\begin{table}[H]
		\centering
		\footnotesize
		\caption{Stability analysis of equilibrium point.}
		\label{st}
		\begin{tabular}{ccccc}
			\hline 
			\multirow{2}{*}{Equilibrium Point} & \multicolumn{2}{c}{Eigenvalues} & \multirow{2}{*}{Stability} & \multirow{2}{*}{Condition} \\
			\cline{2-3}
			& $\lambda(\lambda_{1}, \lambda_{2}, \lambda_{3})$ & Symbol & & \\
			\hline
			$\gamma_{1}(0,0,0)$ & 
			\begin{tabular}{@{}l@{}}
				$\left(C_{S J}-C_{M J}-C_{D J}\right) e^{-\lambda \tau}$ \\
				$C_{I I} e^{-\lambda \tau} $ \\
				$-C_{H J} e^{-\lambda \tau}$
			\end{tabular} &  (*,+, -) & Non-ESS & \\
			
			$\gamma_{2}(1,0,0)$ & 
			\begin{tabular}{@{}l@{}}
				$
				\left(-C_{S J}+C_{M J}+C_{D J}\right) e^{-\lambda \tau} 
				$ \\
				$\left(-C_{SC}+C_{L C}\right) e^{-\lambda \tau} $ \\
				$	C_{I I} e^{-\lambda \tau} $
			\end{tabular} &  (*,*, +) & Non-ESS& \\
			
			$\gamma_{3}(0,1,0)$ & 
			\begin{tabular}{@{}l@{}}
				$\left(C_{S J}-C_{M J}-C_{D J}-C_{H J}-C_{L C}-I_{J}-T_{R J}\right) e^{-\lambda \tau}$ \\
				$C_{H J} e^{-\lambda \tau}$ \\
				$C_{I I} e^{-\lambda \tau}$
			\end{tabular} &  (*,+, +) & Non-ESS & \\
			
			$\gamma_{4}(0,0,1)$ & 
			\begin{tabular}{@{}l@{}}
				$\left(C_{S J}-C_{M J}-C_{D J}\right) e^{-\lambda \tau} $ \\
				$-C_{H J} e^{-\lambda \tau}$ \\
				$-C_{I I} e^{-\lambda \tau}$
			\end{tabular} &  (-,-, -) & ESS & \begin{tabular}{@{}l@{}}
				$C_{SJ}<C_{MJ}+C_{DJ}$
			\end{tabular} \\	
			
			$\gamma_{5}(1,1,0)$ & 
			\begin{tabular}{@{}l@{}}
				$\left(-C_{S J}+C_{M J}+C_{D J}+T_{R J}+C_{H J}+C_{L C}+I_{J}\right) e^{-\lambda \tau} $ \\
				$\left(C_{S C}-C_{L C}\right) e^{-\lambda \tau}$ \\
				$C_{I I} e^{-\lambda \tau}$
			\end{tabular} &  (*,*, +) & Non-ESS & \\

		$\gamma_{6}(1,0,1)$ & 
			\begin{tabular}{@{}l@{}}
				$\left(-C_{S J}+C_{M J}+C_{D J}\right) e^{-\lambda \tau}$ \\
				$\left(C_{L C}-C_{S C}\right) e^{-\lambda \tau}$ \\
				$-C_{I I} e^{-\lambda \tau} $
			\end{tabular} &  (-,-,-) &  ESS& \begin{tabular}{@{}l@{}}
			$C_{S J}>C_{M J}+C_{D J}$ \\
			$ C_{L C}<C_{S C} $ \\
		\end{tabular}\\

				$\gamma_{7}(0,1,1)$ & 
			\begin{tabular}{@{}l@{}}
				$
				\left(C_{S J}-C_{M J}-C_{D J}-C_{H J}-C_{L C}-I_{J}-T_{R J}\right) e^{-\lambda \tau}$ \\
				$C_{H J} e^{-\lambda \tau}$ \\
				$-C_{I I} e^{-\lambda \tau} $
			\end{tabular} &  (*,+,-) & Non-ESS & 	
			\\

			$\gamma_{8}(1,1,1)$ & 
			\begin{tabular}{@{}l@{}}
				$\left(-C_{S J}+C_{M J}+C_{D J}+T_{R J}+C_{H J}+C_{L C}+I_{J}\right) e^{-\lambda \tau}$ \\
				$\left(-C_{L C}+C_{S C}\right) e^{-\lambda \tau}$ \\
				$-C_{I I} e^{-\lambda \tau}$
			\end{tabular} &  (-,-, -) &  ESS &		\begin{tabular}{@{}l@{}}				$C_{L C}>C_{S C}$\\
				$C_{S J}>C_{M J}+C_{D J}+T_{R J}+C_{H J}+C_{L C}+I_{J} $
			\end{tabular} \\
			\bottomrule
		\end{tabular}
		\vspace{6pt}
		\caption*{Note: "-" indicates that the eigenvalue of the Jacobian matrix is negative, "+" means that the eigenvalue is positive, and "*" indicates that the sign of the eigenvalue is uncertain.}
	\end{table}

From Table~\ref{st}, we can see that the equilibria \(\gamma_{1}(0,0,0)\), \(\gamma_{2}(1,0,0)\), \(\gamma_{3}(0,1,0)\), \(\gamma_{5}(1,1,0)\), and \(\gamma_{7}(0,1,1)\) have eigenvalues with positive real parts, hence they are unstable points. The equilibria \(\gamma_{4}(0,0,1)\), \(\gamma_{6}(1,0,1)\), and \(\gamma_{8}(1,1,1)\) have all their eigenvalues with negative real parts under given conditions, making them stable points. We will conduct a preliminary analysis of the situations with these stable points.

\begin{itemize}
	\item For the equilibrium \(\gamma_{4}(0,0,1)\), when the condition \(C_{SJ}<C_{MJ}+C_{DJ}\) is satisfied, it implies that the economic cost for Japan to store nuclear wastewater is less than the combined cost of discharging it to the sea and ocean monitoring. In this scenario, the evolutionary stable point is \(\gamma_{4}(0,0,1)\), where Japan adopts the "stop discharge" strategy, other countries choose the "no sanctions" strategy towards Japan, and IAEA opposes Japan's discharging strategy.
	
	\item For the equilibrium \(\gamma_{6}(1,0,1)\), when the conditions \(C_{SJ}>C_{MJ}+C_{DJ}\) and \(C_{LC}<C_{SC}\) are met simultaneously, it suggests that the economic cost for Japan to store nuclear wastewater exceeds the sum of the costs of discharging it to the sea and ocean monitoring, while other countries' compensations obtained through litigation are less than the extra costs of substituting Japanese products. Under this circumstance, the evolutionary stable point is \(\gamma_{6}(1,0,1)\), where Japan continues to "discharge", other countries adopt a "no sanctions" strategy towards Japan, while IAEA remains "opposing Japan".
	
	\item For the equilibrium \(\gamma_{8}(1,1,1)\), when the conditions \(C_{LC}>C_{SC}\) and \(C_{SJ}>C_{MJ}+C_{DJ}+T_{RJ}+C_{HJ}+C_{LC}+I_{J}\) are both satisfied, it means the cost for Japan to store nuclear wastewater has exceeded the combined costs of sea discharge, ocean monitoring, international reputation loss, reduced export tax, aid from other countries, and compensation from lawsuits to other countries. The extra costs for other countries to develop their own seafood are less than the compensations they get through litigation. In this context, the evolutionary stable point is \(\gamma_{8}(1,1,1)\), where Japan continues to "discharge", other countries enforce a "sanction" strategy against Japan, and the IAEA maintains its "opposing" stance against Japan.
\end{itemize}

In the numerical simulation analysis of Section 6, we will delve into the potential drivers for the formation of evolutionary stable points under the aforementioned scenarios, offering economic and policy interpretations.

\section{Numerical simulation}

In the aforementioned tripartite evolutionary game model, the strategy equilibrium of each game participant is influenced by the strategic choices of the other participants. To more intuitively explore the strategy selection process of the three stable equilibria mentioned above, we analyzed the evolutionary trajectories of stable strategies for the Japanese government, other countries, and the International Atomic Energy Agency. Additionally, simulations were conducted based on different parameter values set for different scenarios.
From Table~\ref{st}, we recognize that the stable points we computed are \(\gamma_{4}(0,0,1)\), \(\gamma_{6}(1,0,1)\), and \(\gamma_{8}(1,1,1)\). To delve deeper into the dynamic characteristics and evolution of these stable points, we will use the replicator dynamics equation and the stability conditions to simulate the evolutionary trajectories of these three stable points. We have set a group of parameter values that satisfy the stability conditions for each stable point. These parameter values will be used to simulate the evolutionary trajectories of these stable points to gain a more accurate understanding of their dynamic behavior. The specific parameter values are shown in Table~\ref{shuzhi}.

	\begin{table}[H]
		\centering
		\caption{Basic parameter settings for three stable equilibrium Points.}
		\label{shuzhi}
		\begin{tabular}{ccccccccccccc}
			\toprule
			Condition&Stable point &$\tau$& $I_{J}$ & $C_{LC}$ & $T_{RJ}$ & $C_{HJ}$ & $C_{DJ}$ & $C_{MJ}$ & $C_{SJ}$ &  $C_{SC}$ & $C_{MC}$ &    $C_{II}$   \\
			\midrule
			Condition 1&$\gamma_{4}(0,0,1)$ &0.01 & 30 & 10 & 10 & 5 & 10 & 10 & 10 & 15 & 10 & 20 \\
			
			Condition 2&$\gamma_{6}(1,0,1)$ &0.01 & 30 & 10 & 10 & 5 & 10 & 10 & 40 & 15 & 10 & 20 \\
			
			Condition 3&$\gamma_{8}(1,1,1)$ &0.01 & 30 & 20 & 10 & 5 & 10 & 10 & 90 & 15 & 10 & 20 \\
			\bottomrule
		\end{tabular}
	\end{table}
	
	\subsection{Evolutionary trajectory under Condition 1}
	\subsubsection{Evolutionary trajectory with $\tau=0$}
	Under Condition 1, where $C_{SJ}<C_{MJ}+C_{DJ}$, this indicates that the cost for Japan to store nuclear wastewater is less than the combined cost of discharging it into the ocean and monitoring the marine environment. In our analysis, we refer to the parameters for Condition 1 as provided in Table~\ref{shuzhi} and set the time delay $\tau=0$.
	
	Firstly, for the Japanese government, when the cost of storing nuclear wastewater is lower than the cost of discharging and monitoring, it can mitigate the negative impacts on marine ecology and its international image by choosing to store the wastewater. This addresses international concerns more effectively. Moreover, considering the long-term impacts of discharging nuclear wastewater and the requirements for sustainable marine development, storing the wastewater can be seen as a delayed decision. This provides more time for the Japanese government to seek more efficient and environmentally friendly solutions, reducing the risk to the marine environment and human health. Thus, the Japanese government is more inclined to adopt a "no discharge" strategy.
	
	For other countries, since the Japanese government adopts the "no discharge" strategy and poses no threats to the marine environment and public health, they might be more inclined to maintain a good cooperative relationship with Japan to further the common goals of ecological safety and environmental protection. Additionally, sanctioning Japan could have adverse effects on their economy and trade. As Japan is a significant economic partner, other countries may rely on trade and investment relations with it. Sanctions could lead to trade disruptions, economic downturns, and employment issues. Thus, to protect their economic interests, these countries might choose a strategy that minimally impacts their national interests, i.e., "not sanctioning" Japan.
	
	For the IAEA, being an organization that coordinates international nuclear cooperation, it would be more concerned about the international community's safety and environmental concerns regarding nuclear wastewater discharge. Even if the Japanese government chooses "no discharge" under $C_{SJ} < C_{MJ} + C_{DJ}$, this does not entirely alleviate public concerns about the impacts of nuclear wastewater and might spark international doubts. In this context, the IAEA would lean more towards encouraging the Japanese government to continue exploring more viable methods for treating nuclear wastewater rather than opting for discharging it. This stance helps maintain the IAEA's international image and responsibilities, ensuring the effective application of science, technology, and environmental protection in the nuclear sector. Hence, the IAEA would choose the "opposition" strategy.
	
\begin{figure}[H]
	\centering
	\begin{minipage}{0.5\textwidth}
		\centering
		\includegraphics[width=0.9\linewidth]{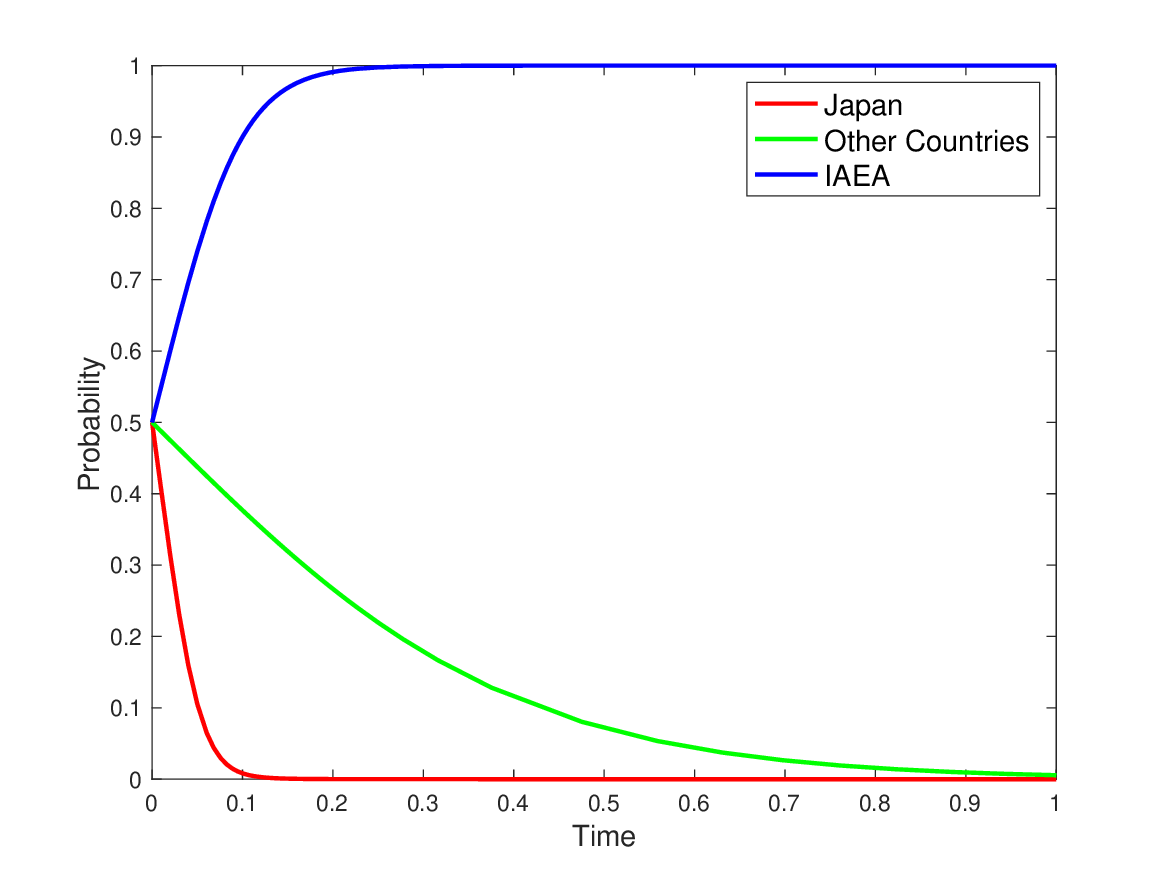}
		\caption{Evolution trajectory over time with $\tau=0$ under \\ Condition 1.}
		\label{c11}
	\end{minipage}%
	\begin{minipage}{0.5\textwidth}
		\centering
		\includegraphics[width=0.9\linewidth]{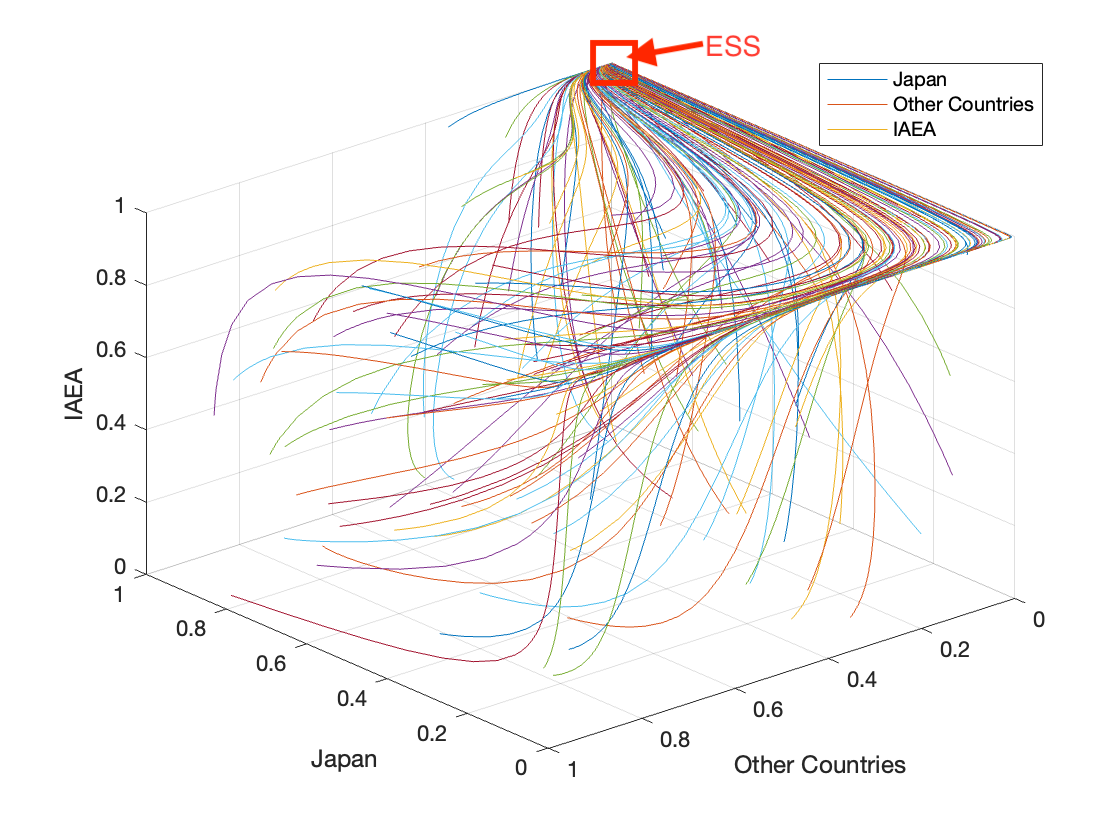}
		\caption{Overall evolution trajectory with $\tau=0$ under \\Condition 1.}
		\label{c12}
	\end{minipage}
\end{figure}

	The results of the numerical simulation further confirm this viewpoint. The simulation shows that under these conditions, the system stabilizes at the point $\gamma_{4}(0,0,1)$. Figure~\ref{c11} displays the evolutionary probabilities of each decision-maker over time, clearly showing Japan's probability of "discharging" approaching 0, other countries' probability of "sanctioning" approaching 0, and the IAEA's "opposition" probability gradually approaching 1. Hence, all the parties tend to stabilize at the point $\gamma_{4}(0,0,1)$ over time. Figure~\ref{c12} displays the different evolutionary trajectories of Japan, other countries, and the IAEA under Condition 1 due to varying initial selection probabilities. Despite being influenced by the parameter ranges, the final decisions of the three parties still tend towards equilibrium, consistent with Figure~\ref{c11}.
	
\subsubsection{Evolution trajectory with different time delays $\tau$}

\begin{figure}[H]
	\centering
	\begin{subfigure}{0.33\textwidth}
		\centering
		\includegraphics[width=1\linewidth]{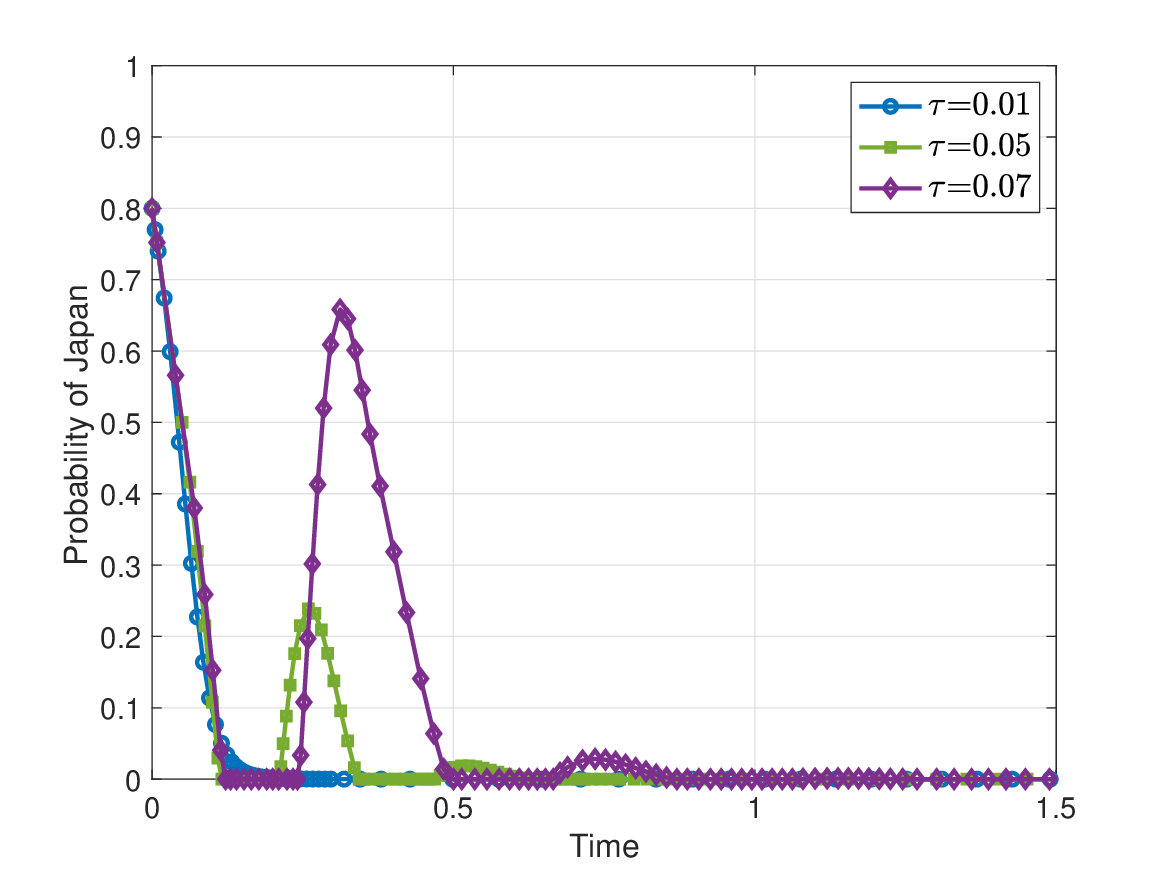}
		\caption{Impact of time delay $\tau$ on Japan.}
		\label{tau1}
	\end{subfigure}%
	\begin{subfigure}{0.33\textwidth}
		\centering
		\includegraphics[width=1\linewidth]{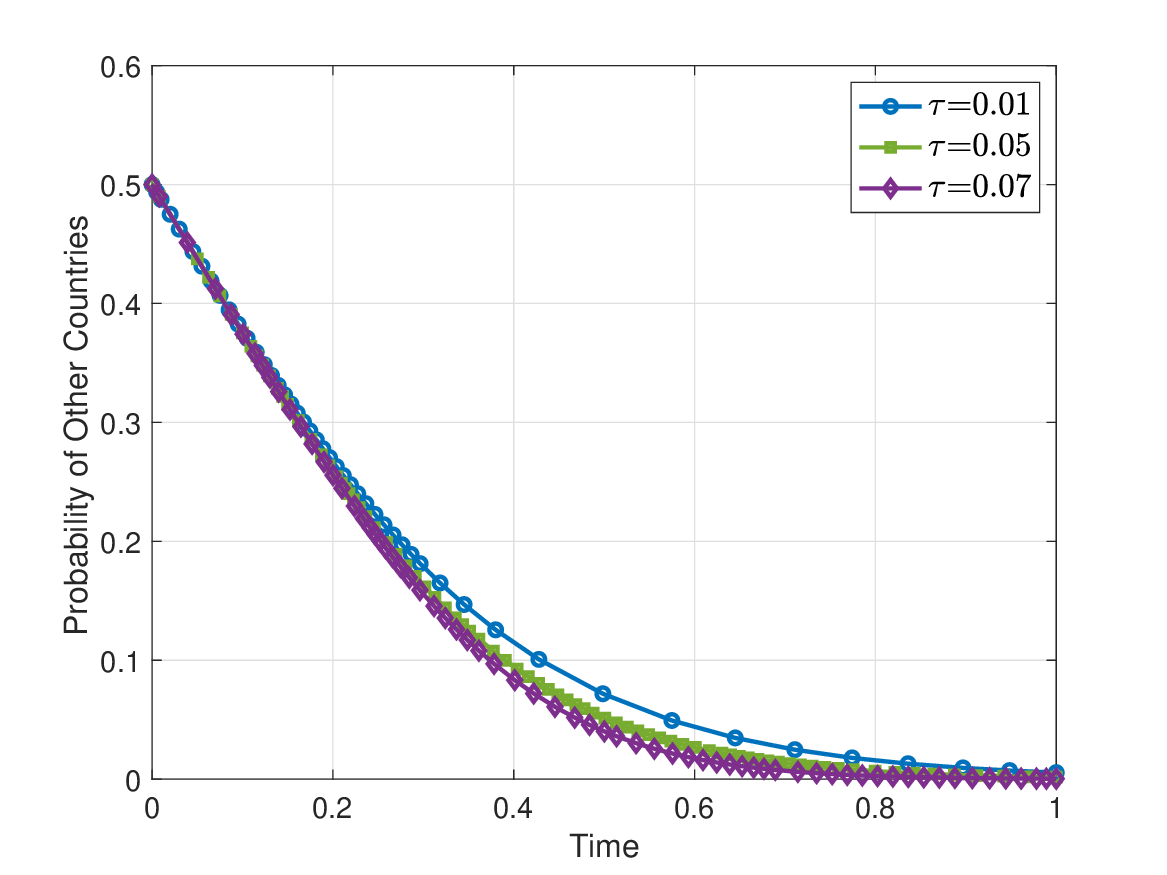}
		\caption{Impact of time delay $\tau$ on other countries.}
		\label{tau2}
	\end{subfigure}%
	\begin{subfigure}{0.33\textwidth}
		\centering
		\includegraphics[width=1\linewidth]{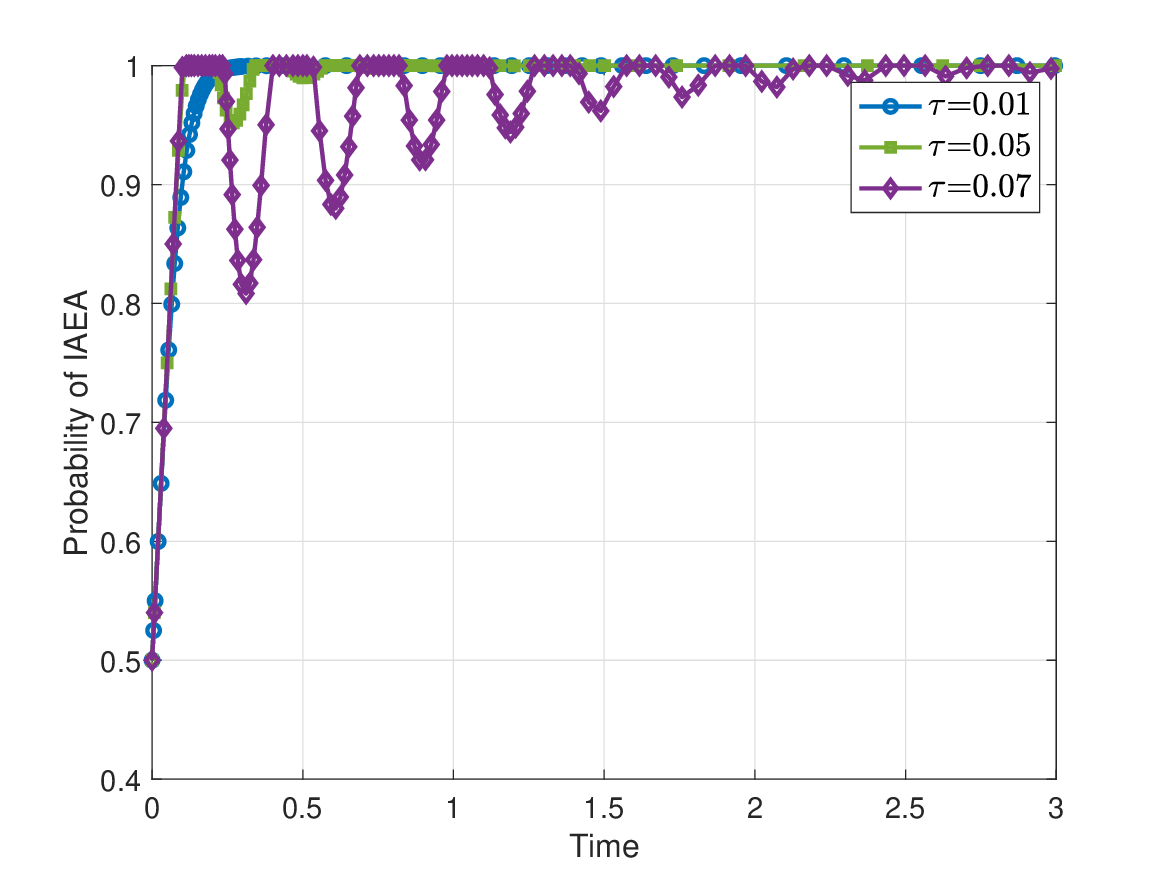}
		\caption{Impact of time delay $\tau$ on IAEA.}
		\label{tau3}
	\end{subfigure}
	\caption{Impact of different time delays $\tau$ on evolution trajectories.}
	\label{tau123}
\end{figure}

Under Condition 1, with [0.8,0.5,0.5] as the initial values and time delays $\tau$ set to 0.01, 0.05, and 0.07 respectively, the numerical simulation reveals the evolution trajectories for the three game players as shown in Figure~\ref{tau123}.

From Figure~\ref{tau123}, it can be seen that although Japan, other countries, and the IAEA are all influenced by the time delay, the time to reach equilibrium is delayed. The longer the delay, the longer the system takes to stabilize. However, the final equilibrium decision remains consistent with the analysis in section 6.1.1. Under this condition:

The Japanese government will adjust its decisions more frequently during the game. As $\tau$ increases, the time taken for Japan's evolution to reach a stable state also increases. This indicates that the Japanese government, considering caution, will pay more attention to risk management regarding the discharge of nuclear wastewater, especially when the storage costs are relatively low. Thus, the Japanese government prefers the "no-discharge" strategy, as shown in Figure~\ref{tau1}.

The time delay has a minimal impact on the evolution trajectory of other countries, meaning they will choose options with the least national interest impact, i.e., "not sanctioning" Japan, as seen in Figure~\ref{tau2}.

For the IAEA, as time delay increases, the frequency of strategy adjustments also increases. The IAEA, placing importance on future marine sustainable development and international public opinion, needs to repeatedly review and base its decisions on the conditions of Japan's nuclear wastewater discharge. The IAEA ultimately chooses to "oppose" Japan's discharge strategy, as illustrated in Figure~\ref{tau3}.

\subsection{Evolution trajectory under Condition 2}
\subsubsection{Evolution trajectory with $\tau=0$}
When the system satisfies Condition 2, that is, $C_{S J}>C_{M J}+C_{D J}$ and $ C_{L C}<C_{S C} $. This condition implies that the cost for Japan to store nuclear wastewater exceeds the sum of the cost of discharging into the sea and the cost of marine monitoring. Additionally, the compensation from a lawsuit for other countries is less than the extra cost of substituting Japanese products. We set time delay $\tau=0$.

When the Japanese government weighs the strategies of nuclear wastewater storage versus discharging, it considers the costs and benefits of wastewater treatment. Storing nuclear wastewater requires a significant amount of funds and resources, and long-term storage poses challenges. In contrast, discharging nuclear wastewater into the sea can avoid the high costs associated with storage and treatment. Hence, if the cost of storing nuclear wastewater exceeds the combined cost of discharging and marine monitoring, the Japanese government might find discharging to be more economically beneficial, opting for the "discharge" strategy.

If the lawsuit compensation received by other countries is less than the additional cost required to substitute Japanese products, even in the face of controversies over wastewater discharge, other countries might prioritize economic benefits and choose "no sanctions". 
(i) When other countries consider substituting Japanese products, they will evaluate potential technological transitions, supply chain adjustments, and the reselection of business partners, all of which may incur additional costs. If the substitution cost is high, they may prefer to continue purchasing Japanese products.
(ii) If other countries have a close economic relationship with Japan and Japan is a critical product supplier, sanctions could lead to supply interruptions or create supply-demand gaps, severely affecting their economy. When other countries find it challenging to quickly locate alternative sources, this supply-demand relationship and mutual dependency become significant factors in their decision to "not sanction Japan". Thus, under this backdrop, other countries ultimately choose "no sanctions".

Protecting the environment, ecology, and human health are important concerns shared by the international community. International standards and collaborative frameworks also emphasize sustainable development and nuclear safety. Given the focus on risks and environmental impacts of nuclear wastewater discharge, the IAEA, in consensus with the international community, emphasizes protecting marine ecosystems and human life, ultimately choosing the "oppose discharge" strategy. This choice reinforces the IAEA's image as an authoritative institution, attentive to global environmental issues, promoting sustainable development, and strengthening international cooperation. It further solidifies the IAEA's leadership position in the nuclear energy sector and boosts its public credibility.

\begin{figure}[H]
		\centering
		\begin{minipage}{0.5\textwidth}
			\centering
			\includegraphics[width=0.9\linewidth]{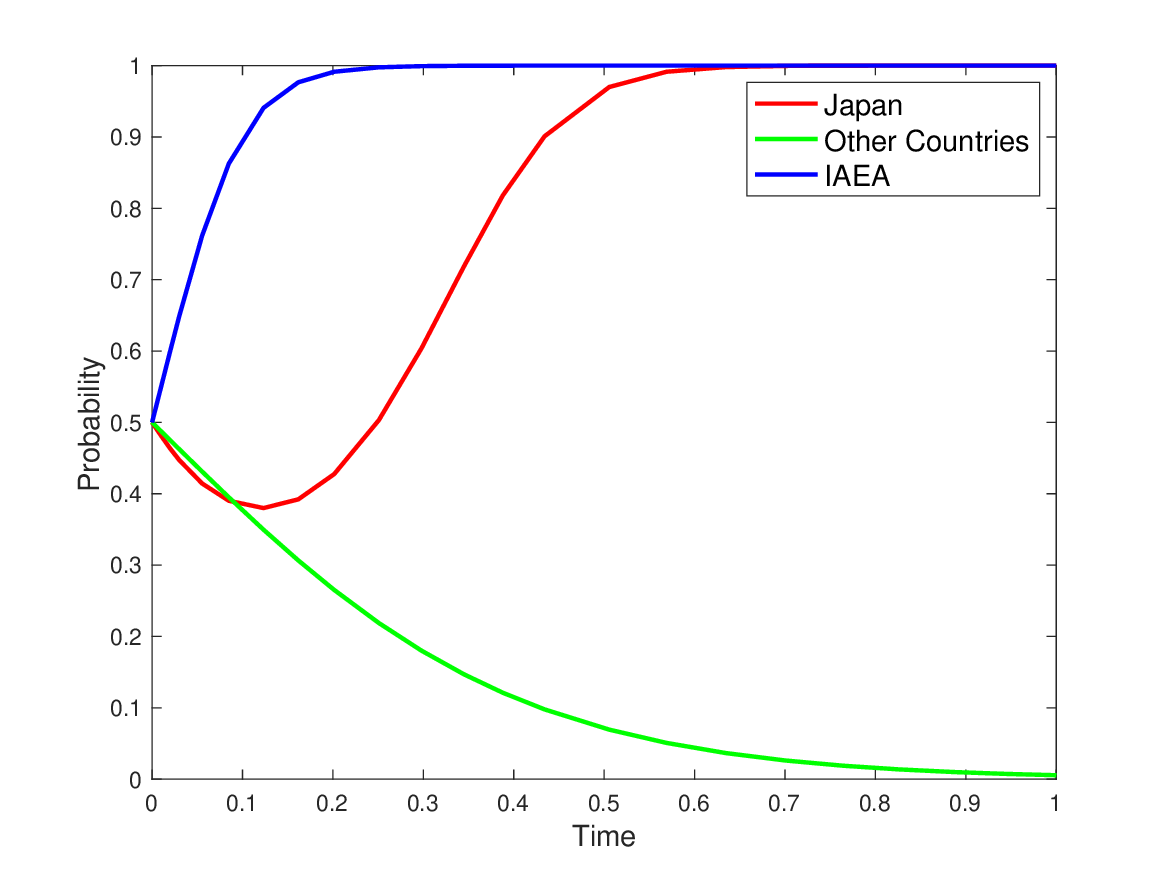}
			\caption{Evolution trajectory over time under Condition 2 with $\tau=0$.}
			\label{c21}
		\end{minipage}%
		\begin{minipage}{0.5\textwidth}
			\centering
			\includegraphics[width=0.9\linewidth]{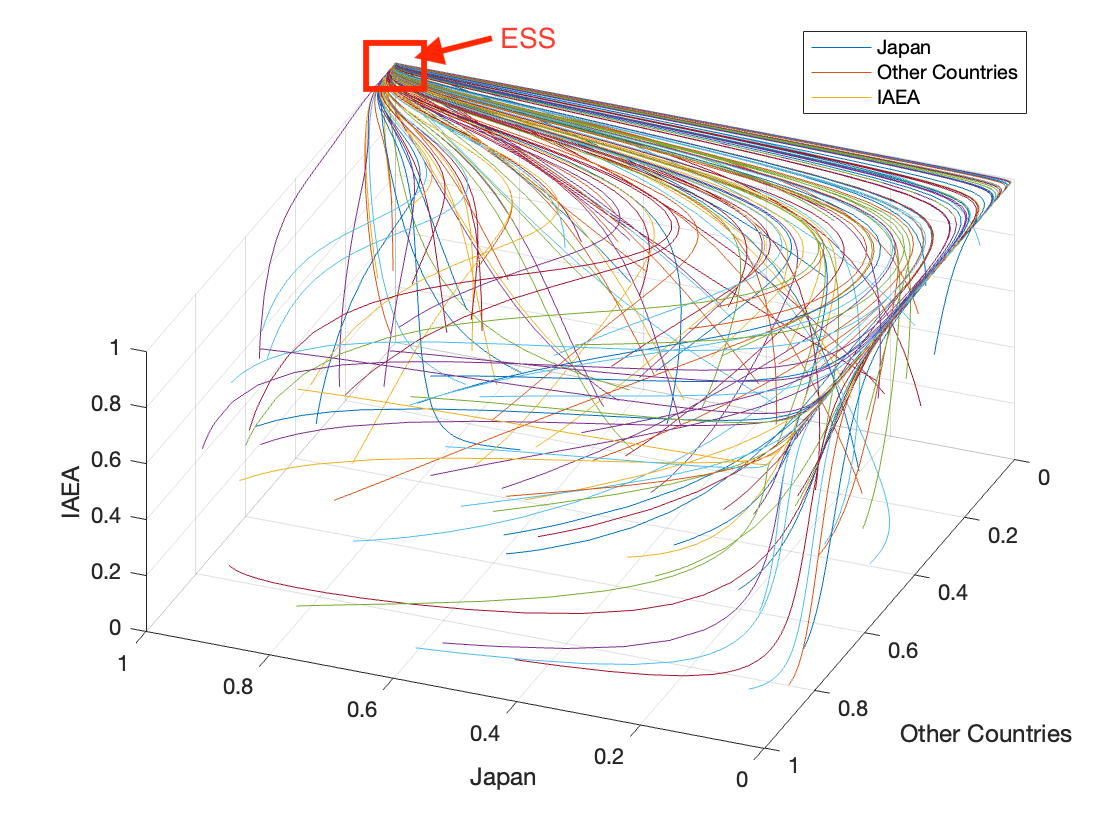}
			\caption{Overall evolutionary trajectory under Condition 2 with $\tau=0$.}
			\label{c22}
		\end{minipage}
	\end{figure}

	Numerical simulation results further confirm this view. The simulation results show that, under these conditions, the system's equilibrium point is $\gamma_{6}(1,0,1)$. Figure~\ref{c21} displays how the evolution probabilities of each decision-maker change over time. It clearly shows the probability of Japan choosing "discharge" gradually approaching 1, the probability of other countries choosing "sanctions" tending to 0, and the IAEA's "opposition" probability gradually approaching 1. This means each entity's decisions are gradually moving towards the equilibrium point $\gamma_{6}(1,0,1)$. Figure ~\ref{c22}, on the other hand, displays the different evolutionary trajectories of Japan, other countries, and the IAEA under Condition 2 due to different initial choice probabilities. Influenced by parameter constraints, the final decisions of all three parties will still tend to balance, which corresponds with Figure~\ref{c21}.
	
	\subsubsection{Evolutionary trajectories under different time delays $\tau$}
	
\begin{figure}[H]
		\centering
		\begin{subfigure}{0.33\textwidth}
			\centering
			\includegraphics[width=1\linewidth]{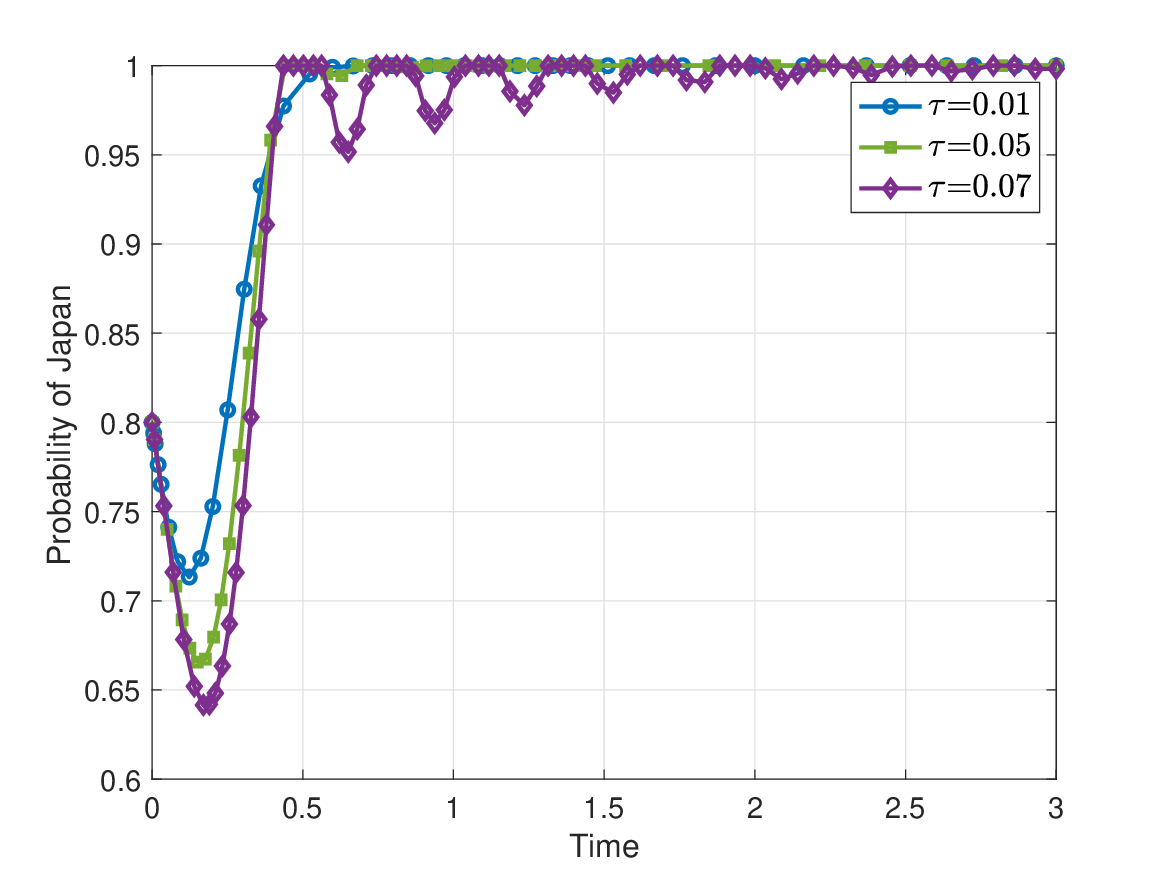}
			\caption{Impact of time delay $\tau$ on Japan.}	
			\label{tau4}
		\end{subfigure}%
		\begin{subfigure}{0.33\textwidth}
			\centering
			\includegraphics[width=1\linewidth]{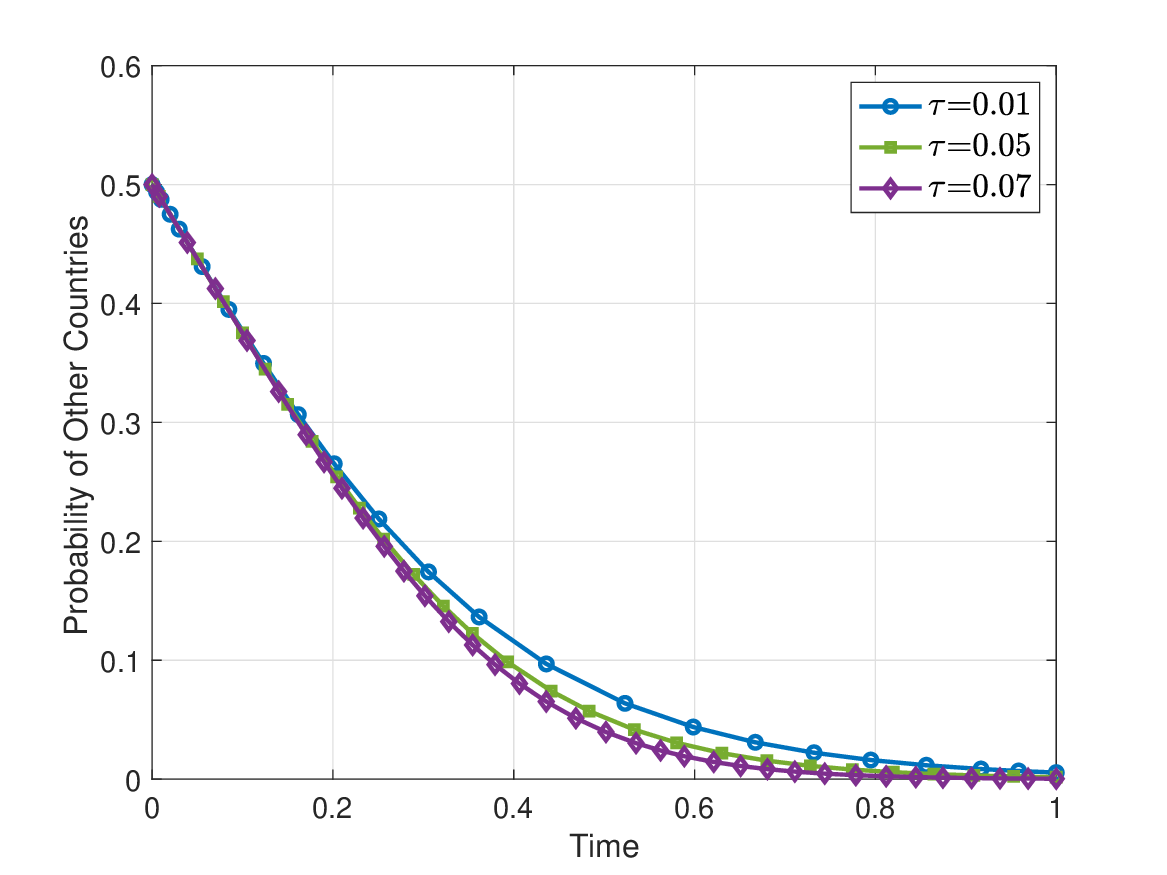}
			\caption{Impact of time delay $\tau$ on other countries.}
			\label{tau5}
		\end{subfigure}%
		\begin{subfigure}{0.33\textwidth}
			\centering
			\includegraphics[width=1\linewidth]{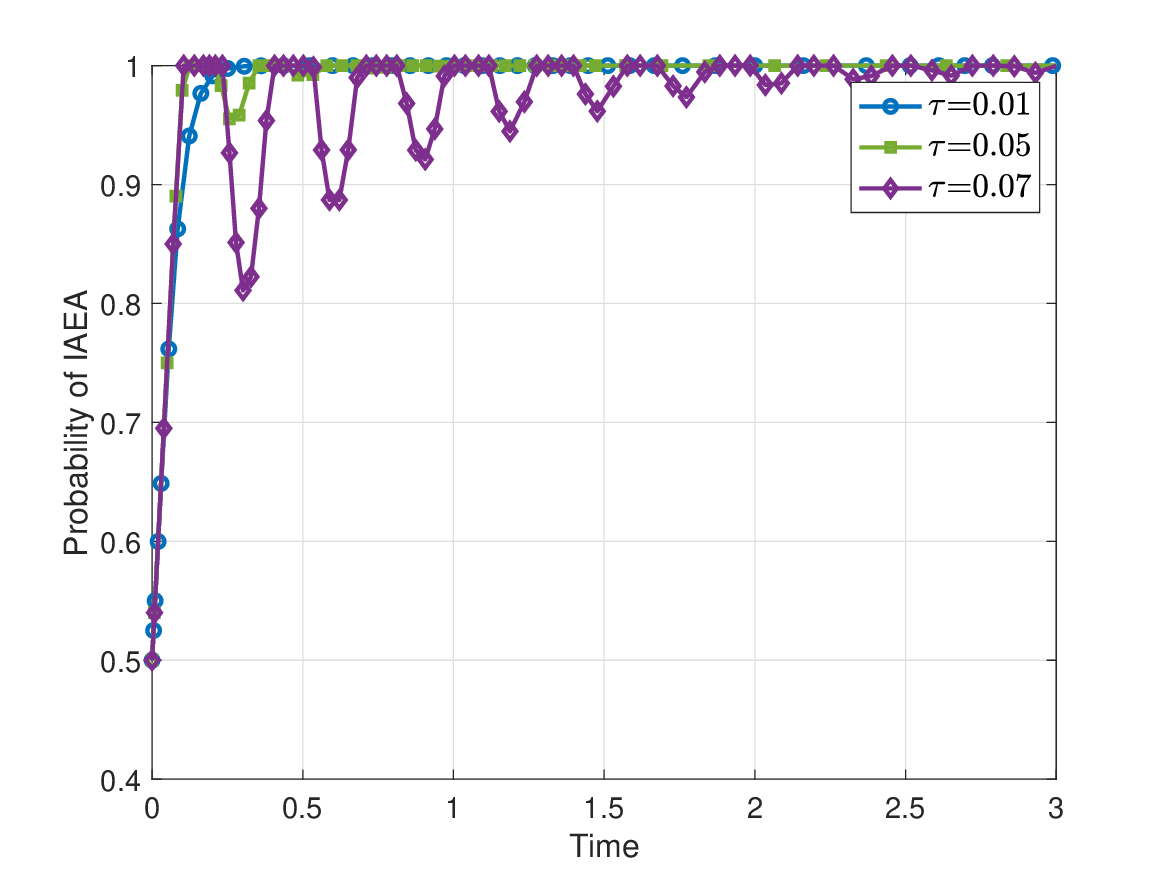}
			\caption{Impact of time delay $\tau$ on IAEA.}
			\label{tau6}
		\end{subfigure}
		\caption{Effects of different time delays $\tau$ on evolutionary trajectories.}
		\label{tau456}
	\end{figure}

	Under Condition 2, starting from the initial values [0.8,0.5,0.5] and setting the time delay $\tau$ to 0.01, 0.05, and 0.07 respectively, the numerical simulations show the effects of different time delays $\tau$ on the evolutionary trajectories of the three game participants as illustrated in Figure~\ref{tau456}.
	
	From Figure~\ref{tau456}, it can be seen that even though Japan, other countries, and the IAEA are all affected by the time delay during the game process, the final decision is consistent with the analysis in 6.2.1. Under Condition 2, as the Japanese government continues to adjust its decision scheme throughout the game process, the number of adjustments increases. The larger the time delay $\tau$, the longer it takes for the Japanese government to evolve to its final stable state. This indicates that, after weighing the options, the Japanese government still believes that discharging the nuclear wastewater is more economically beneficial. Thus, the Japanese government ultimately chooses the "discharge" strategy, as shown in Figure~\ref{tau4}.
	
	The time delay has a minor effect on other countries' strategy adjustments. In this scenario, other countries tend to stick to their original decisions rather than adjust rapidly. This might be because their decision-making process involves various factors, including the litigation compensation they receive, the additional costs of replacing Japanese products, economic trade, and international relations. These factors still play a role in the presence of time delay, leading other countries to lean towards maintaining their initial strategies, as depicted in Figure~\ref{tau5}.
	
	Due to the time delay, for the IAEA, the number of strategy adjustments increases with the increase of delay. Also, the time it takes for the IAEA to evolve to its final stable state extends with the time delay $\tau$. This indicates that the IAEA values the sustainable development of the future ocean and international public opinion. Their decision-making requires a comprehensive consideration and repeated validation. Eventually, the IAEA opts for the strategy of "opposing Japan's discharge", as shown in Figure~\ref{tau6}.
	
\subsection{Evolutionary trajectory under Condition 3}
\subsubsection{Evolutionary trajectory with $\tau=0$}

When the system meets Condition 3, i.e., $C_{L C} > C_{S C}$ and $C_{S J} > C_{M J} + C_{D J} + T_{R J} + C_{H J} + C_{L C} + I_{J}$, we refer to the parameters of Condition 3 in Table \ref{shuzhi} for analysis. The cost for Japan to store nuclear wastewater exceeds the costs of ocean discharge, marine monitoring, international image damage, reduced export taxes, and the aid obtained from other countries, as well as litigation compensation to other countries. The additional cost for other countries to develop their marine products is less than the litigation compensation they receive. We set the time delay $\tau=0$.

The cost for Japan to store nuclear wastewater has already surpassed the costs of ocean discharge, marine monitoring, damage to international image, reduced export tax revenue, and both the aid from other countries and the compensation paid to them for litigation. By discharging the nuclear wastewater into the ocean, the government can avoid the expensive costs required for storing and processing it. Moreover, as the nuclear wastewater issue becomes internationally exposed, Japan's image will face negative reviews, which in turn leads to reduced export trade and tax revenue. The Japanese government will weigh these losses against the benefits of discharging the wastewater, ultimately choosing the option with the least negative impact on the nation's overall interests. This means the Japanese government is more inclined to choose the "discharge" strategy to deal with the nuclear wastewater.

For other countries, they tend to choose to sanction Japan's discharge for the following reasons: (i) There may be additional costs involved in the process of importing or producing alternatives, increasing the extra expenses of other countries. (ii) Other countries might evaluate the potential benefits they could reap by importing these alternatives. Therefore, in order to minimally affect the nation's overall interests, other countries will choose to "sanction" Japan's discharge.

As for the International Atomic Energy Agency (IAEA), when the Japanese government decides to discharge nuclear wastewater into the sea, the IAEA adopts an opposing strategy. Firstly, the IAEA is committed to the safe and sustainable development of nuclear energy, prioritizing environmental protection and nuclear radiation risk management. Discharging nuclear wastewater poses potential risks to marine ecosystems and biodiversity, which conflicts with the IAEA's nuclear safety objectives. Secondly, as an institution coordinating international nuclear cooperation, the international public is widely concerned about the safety and environmental impacts of nuclear wastewater discharge. The IAEA is more likely to adopt a cautious stance on the discharge, based on maintaining its international image and upholding its scientific and technical authority. Thus, the IAEA tends to "oppose Japan's discharge".

\begin{figure}[H]
	\centering
	\begin{minipage}{0.5\textwidth}
		\centering
		\includegraphics[width=0.9\linewidth]{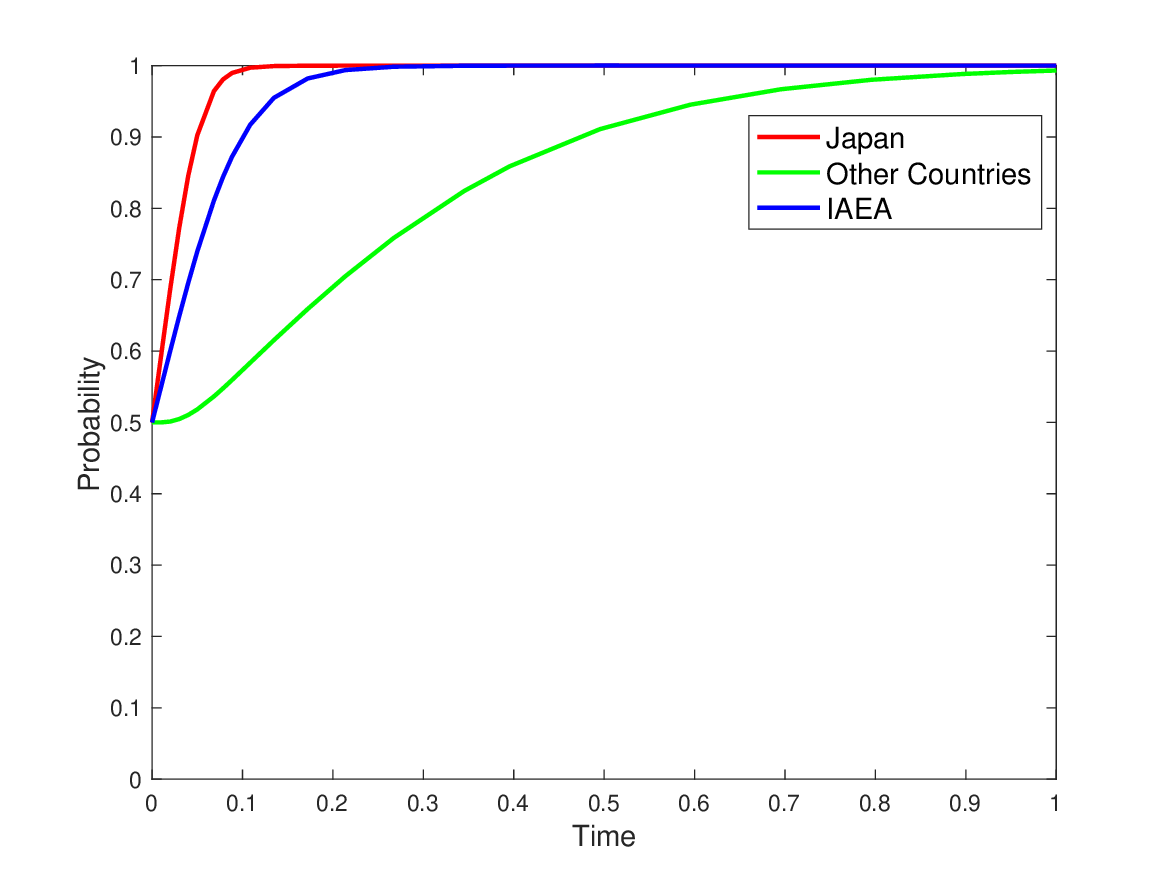}
		\caption{Evolutionary trajectory over time under \\Condition 3 with $\tau=0$.}
		\label{c31}
	\end{minipage}%
	\begin{minipage}{0.5\textwidth}
		\centering
		\includegraphics[width=0.9\linewidth]{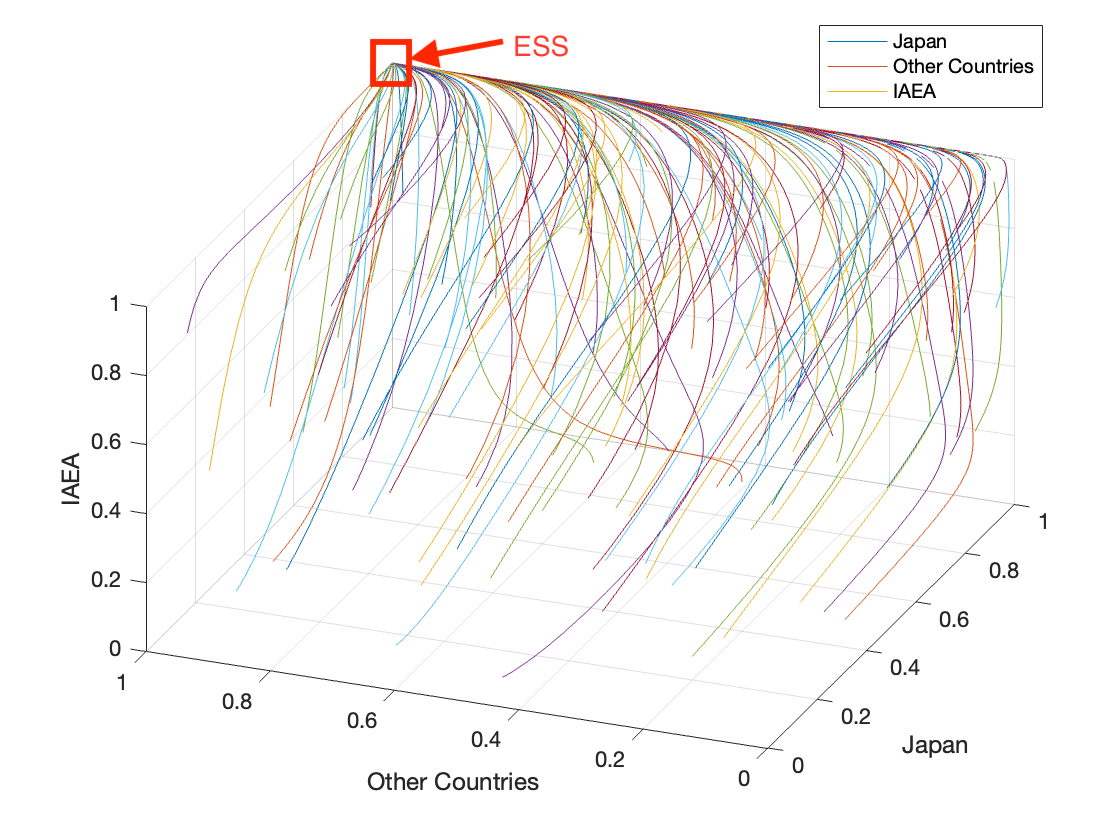}
		\caption{Overall evolutionary trajectory under Condition 3 with $\tau=0$.}
		\label{c32}
	\end{minipage}
\end{figure}

Numerical simulations further confirm this perspective. The simulation results show that under these conditions, the system's stable point is $\gamma_{8}(1,1,1)$. Figure~\ref{c31} depicts the evolution of probabilities of each decision-making entity over time. It is clear that the probability of Japan choosing "discharge" tends to 1, the probability of other countries choosing "sanctions" tends to 1, and the probability of the IAEA choosing "opposition" tends to 1. This indicates that over time, each entity tends to the stable point $\gamma_{8}(1,1,1)$. Figure~\ref{c32} displays the different evolutionary trajectories under Condition 3 for Japan, other countries, and the IAEA due to differences in initial choice probabilities. Despite being influenced by the parameter constraints, the final decisions of all three parties still tend towards equilibrium, corresponding to Figure~\ref{c31}.

\subsubsection{Evolutionary trajectories under different time delays $\tau$}

\begin{figure}[H]
	\centering
	\begin{subfigure}{0.33\textwidth}
		\centering
		\includegraphics[width=1\linewidth]{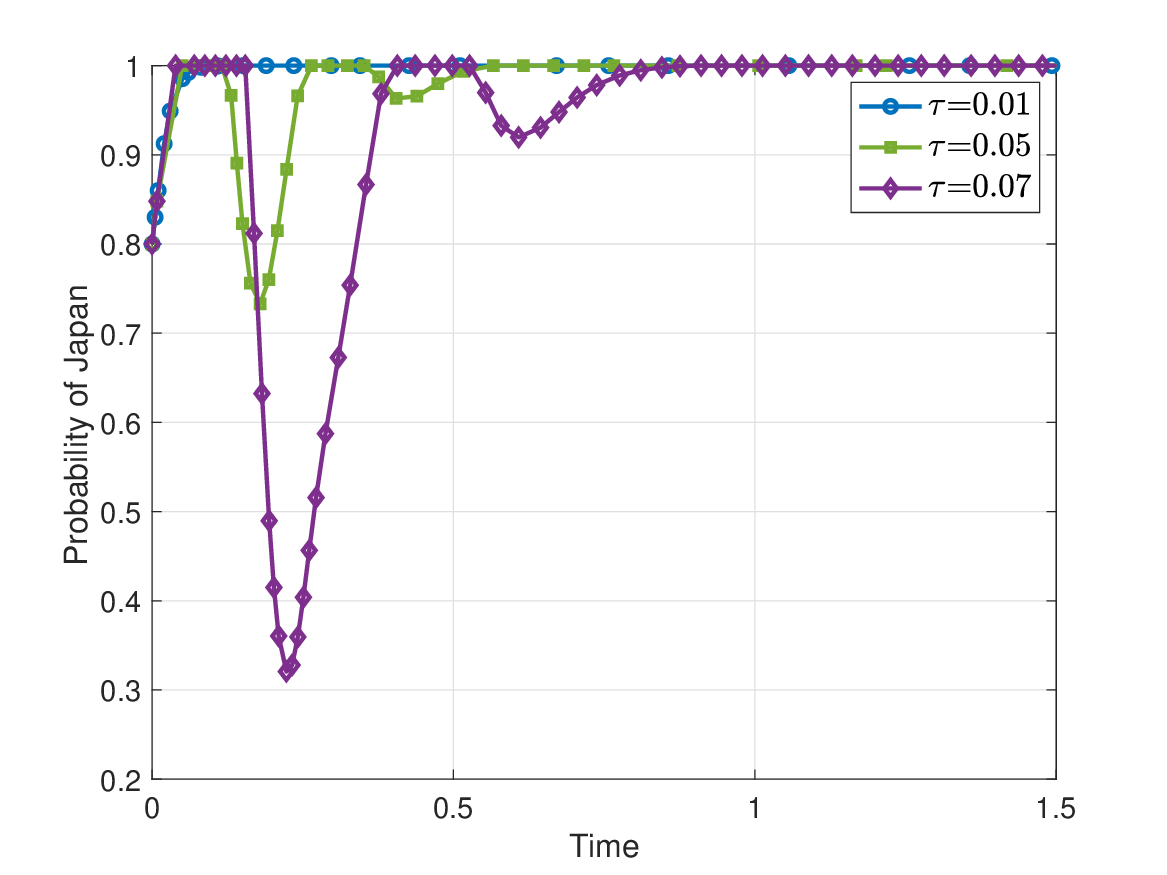}
		\caption{Effect of time delay $\tau$ on Japan.}
		\label{tau7}
	\end{subfigure}%
	\begin{subfigure}{0.33\textwidth}
		\centering
		\includegraphics[width=1\linewidth]{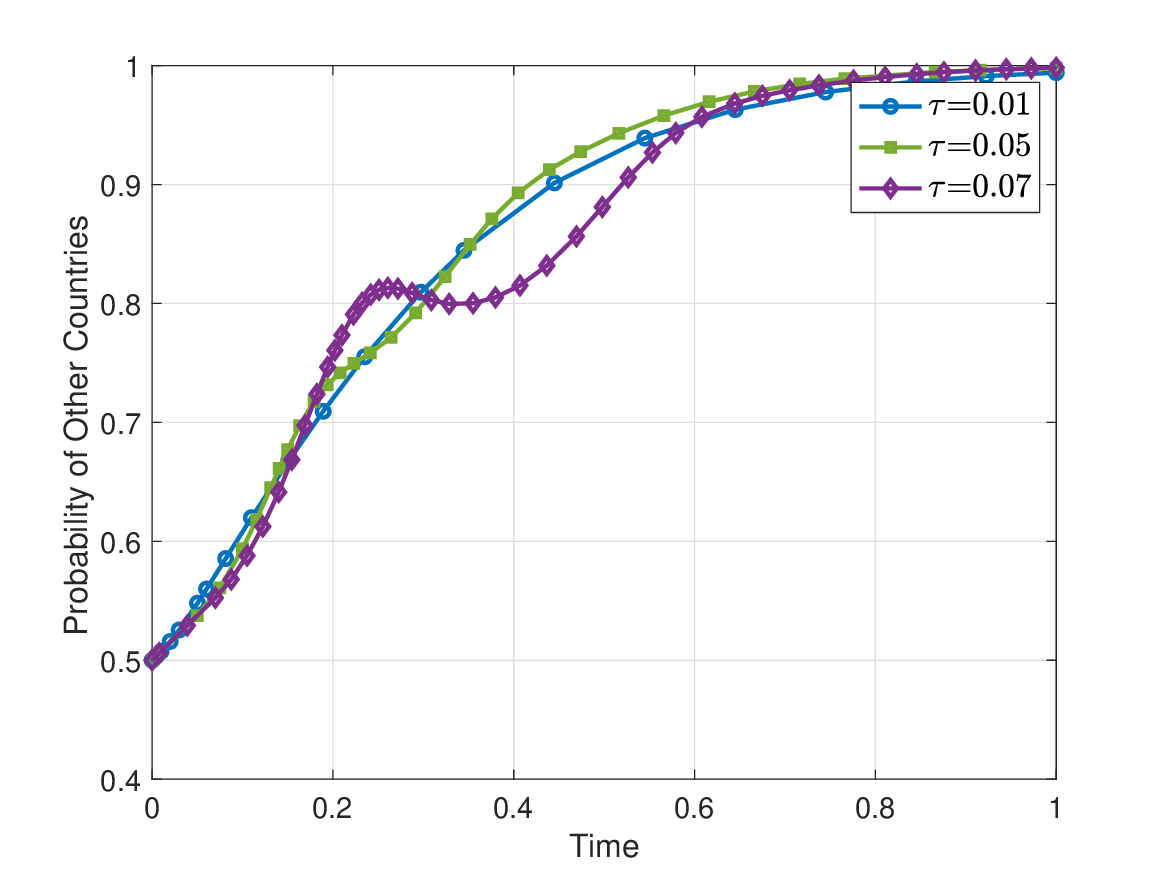}
		\caption{Effect of time delay $\tau$ on other countries.}
		\label{tau8}
	\end{subfigure}%
	\begin{subfigure}{0.33\textwidth}
		\centering
		\includegraphics[width=1\linewidth]{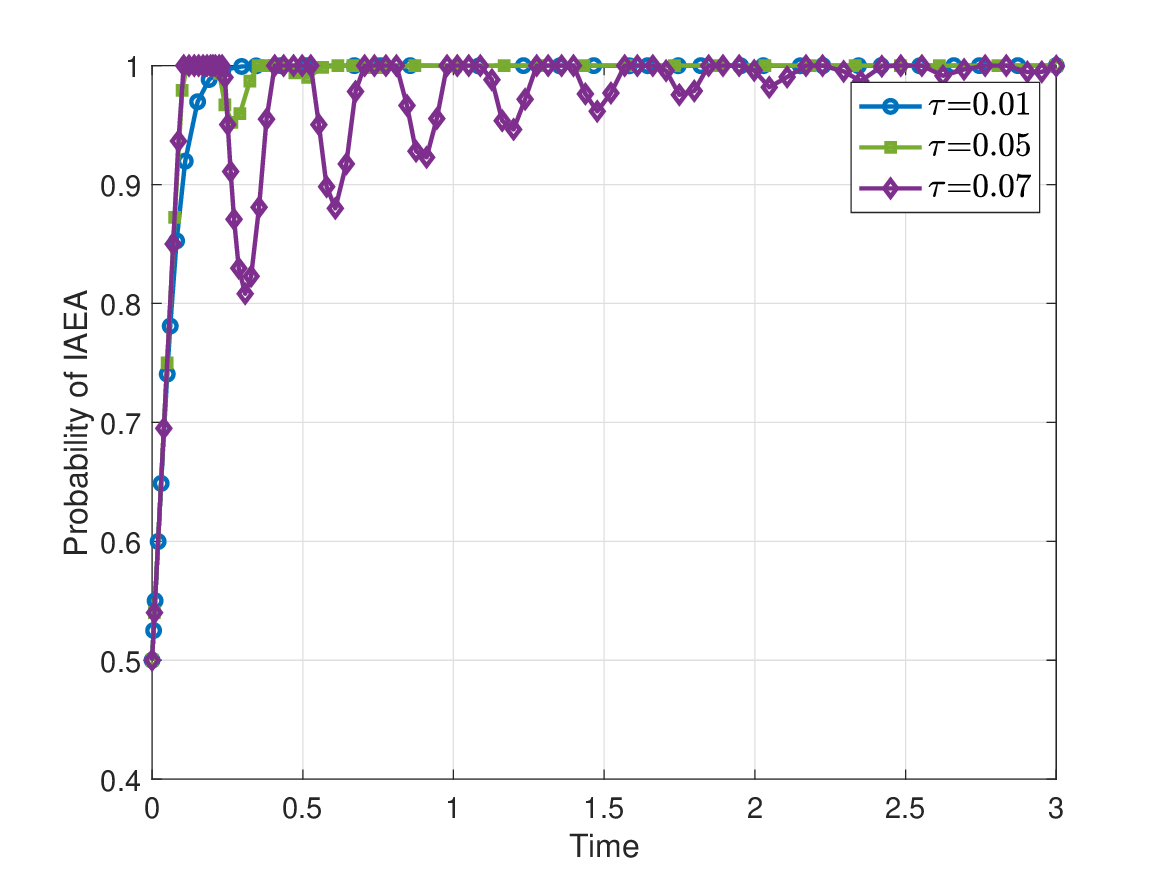}
		\caption{Effect of time delay $\tau$ on IAEA}
		\label{tau9}
	\end{subfigure}
	\caption{Impact of different time delays $\tau$ on evolutionary trajectories.}
	\label{tau789}
\end{figure}

Under Condition 3, using [0.8,0.5,0.5] as the initial values and setting the time delays $\tau$ to 0.01, 0.05, and 0.07 respectively, the numerical simulations provide the evolutionary trajectories for the three game participants as shown in Figure~\ref{tau789}.

As can be observed in Figure~\ref{tau789}, although Japan, other countries, and the IAEA are all influenced by the time delay during the game process, the final decision-making outcome is consistent with the analysis in section 6.3.1. Under Condition 3, due to the time delay, the frequency of adjustments in Japan's decisions during the game process increases. The larger the time delay $\tau$, the longer it takes for Japan to evolve to its final stable state, as depicted in Figure~\ref{tau7}. This indicates that, even though Japan considers the potential negative impact on its national image from international exposure, leading to a decrease in export transactions and tax revenue, the Japanese government still leans towards the "discharge" strategy due to the high storage costs.

Time delay causes more evident changes in the evolutionary trajectories of other countries. As the time delay τ increases, these countries give more consideration to the potential additional costs of introducing or producing alternatives. Thus, after weighing their options to minimally impact national interests, these countries will choose to "sanction" Japan's discharge, as illustrated in Figure~\ref{tau8}.

For the IAEA, under the influence of time delay, the number of strategy adjustments increases with the delay. Moreover, as the time delay τ grows, the time it takes for the IAEA to evolve to its final stable state also elongates. This indicates that the IAEA gives more importance to the sustainable development of the oceans in the future and the opinion of the international public. Ultimately, the IAEA opts for the "oppose" Japan's discharge strategy, as seen in Figure~\ref{tau9}.

\subsection{Impact of different parameter combinations on evolutionary trajectories}
Our focus is on Japan's final strategy choice to stop discharging, i.e., Condition 1. To delve deeper into the potential effects of certain parameters such as Japan's discharge cost (\(C_{DJ}\)), the cost for Japan to store nuclear wastewater (\(C_{SJ}\)), the litigation compensation from other countries (\(C_{LC}\)), assistance Japan receives from other countries (\(C_{HJ}\)), and the reduction in export tax revenue for Japan due to discharging ($T_{RJ}$) on the stable strategies that might evolve among the three parties, a series of numerical simulations were conducted under the conditions described in Table~\ref{biao5}.

	\begin{table}[h]
		\centering
		\caption{Basic parameter settings.}
		\begin{tabular}{cccccccccccc}
	\toprule
	Condition&$\tau$& $I_{J}$ & $C_{LC}$ & $T_{RJ}$ & $C_{HJ}$ & $C_{DJ}$ & $C_{MJ}$ & $C_{SJ}$ &  $C_{SC}$ & $C_{MC}$ &    $C_{II}$   \\
	\midrule
	Condition 4 &0.01 & 30 & 10 & 10 & 5 &  [1,30]  & 10 & [1,30] & 15 & 10 & 20 \\
		\midrule
	Condition 5 &0.01 & 30 & 10 & [1,30]  & 5 & 10 & 10 & [1,30] & 15  & 10 & 20 \\
		\midrule
	Condition 6 &0.01 & 30 & [1,18] & 10 & [1,30] & 10 & 10 & 80 & 15 & 10 & 20 \\
			\bottomrule			\label{biao5}
		\end{tabular}	
	\end{table}

\subsubsection{Analysis of the impact of Japan's nuclear wastewater storage cost $C_{SJ}$ and discharge cost $C_{DJ}$ on decision-making attitude}

\begin{figure}[H]
	\centering
	\includegraphics[width=0.6\linewidth]{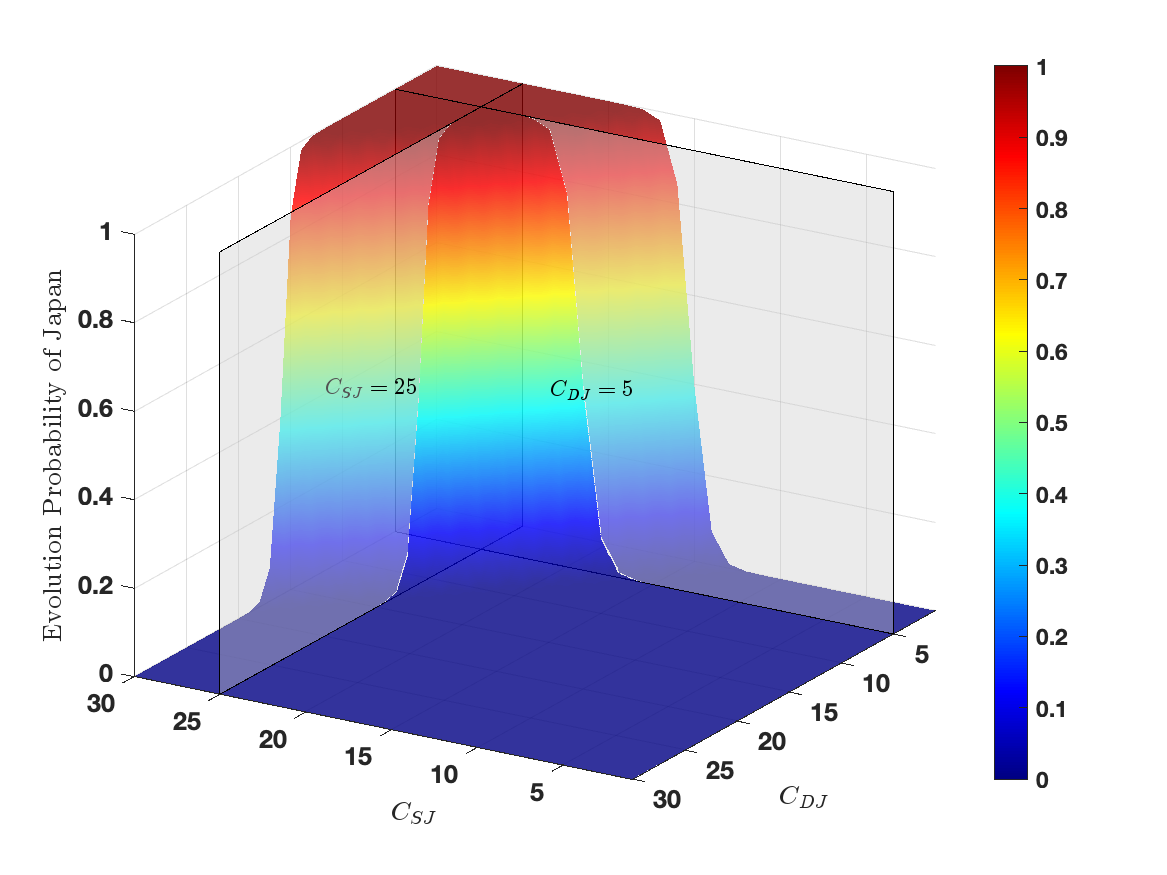}
	\caption{Impact of $C_{SJ}$ and $C_{DJ}$ on Japan's evolutionary decision-making.}
	\label{duo1}
\end{figure}

The handling of nuclear wastewater, whether storage or discharge, involves substantial economic and environmental costs. Against this backdrop, the storage cost of nuclear wastewater $C_{SJ}$ and the discharge cost $C_{DJ}$ become pivotal factors determining the stance and decisions of the Japanese government. To gain a more comprehensive understanding of how these two cost factors influence Japan's decision-making, we referred to the initial parameter settings of Condition 4 in Table~\ref{biao5} and drew a three-dimensional surface plot to visually showcase how these two cost variables shape Japan's policy attitude, as depicted in Figure~\ref{duo1}.

Figure~\ref{duo1} lucidly displays the effects of $C_{SJ}$ and $C_{DJ}$ on Japan's decisions: as $C_{SJ}$ increases, meaning the cost of storing nuclear wastewater rises, Japan's policy inclination gradually shifts from "no discharge" of nuclear wastewater to "discharge". This is because, as storage costs mount, discharging becomes a more economically viable option. Conversely, when $C_{DJ}$ rises, signifying an increase in the costs of discharging nuclear wastewater, Japan's policy inclination moves from "discharge" to "no discharge". This indicates that when the economic and environmental costs of discharging become too steep, the Japanese government would reconsider its decision. Specifically, under the aforementioned parameter settings, when $C_{DJ}$ exceeds 20 or $C_{SJ}$ is less than 10, the Japanese government consistently opts for the "no discharge" strategy.

\begin{figure}[H]
	\centering
	\begin{subfigure}{0.45\textwidth}
		\centering
		\includegraphics[width=1\linewidth]{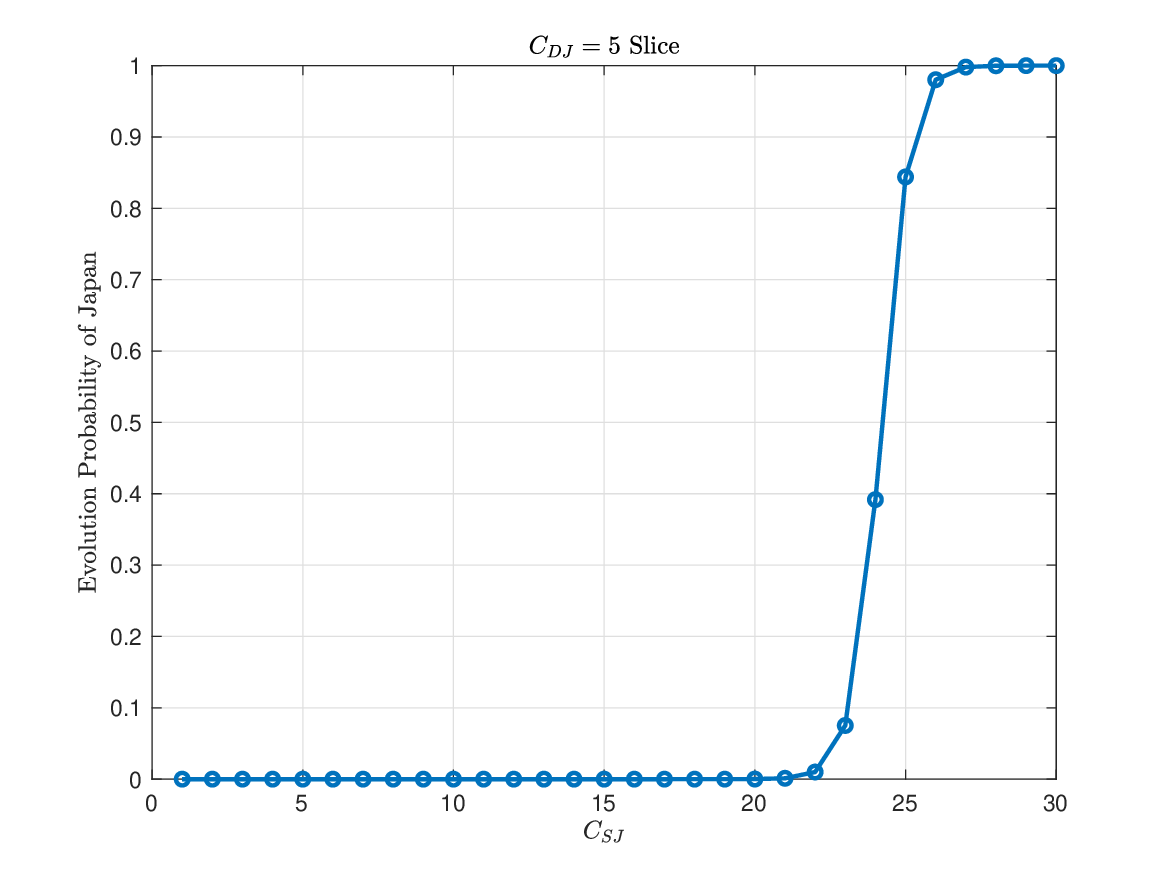}
		\caption{Impact of $C_{SJ}$ on Japan's Evolutionary Decision-making.}
		\label{duo3}
	\end{subfigure}%
	\begin{subfigure}{0.45\textwidth}
		\centering
		\includegraphics[width=1\linewidth]{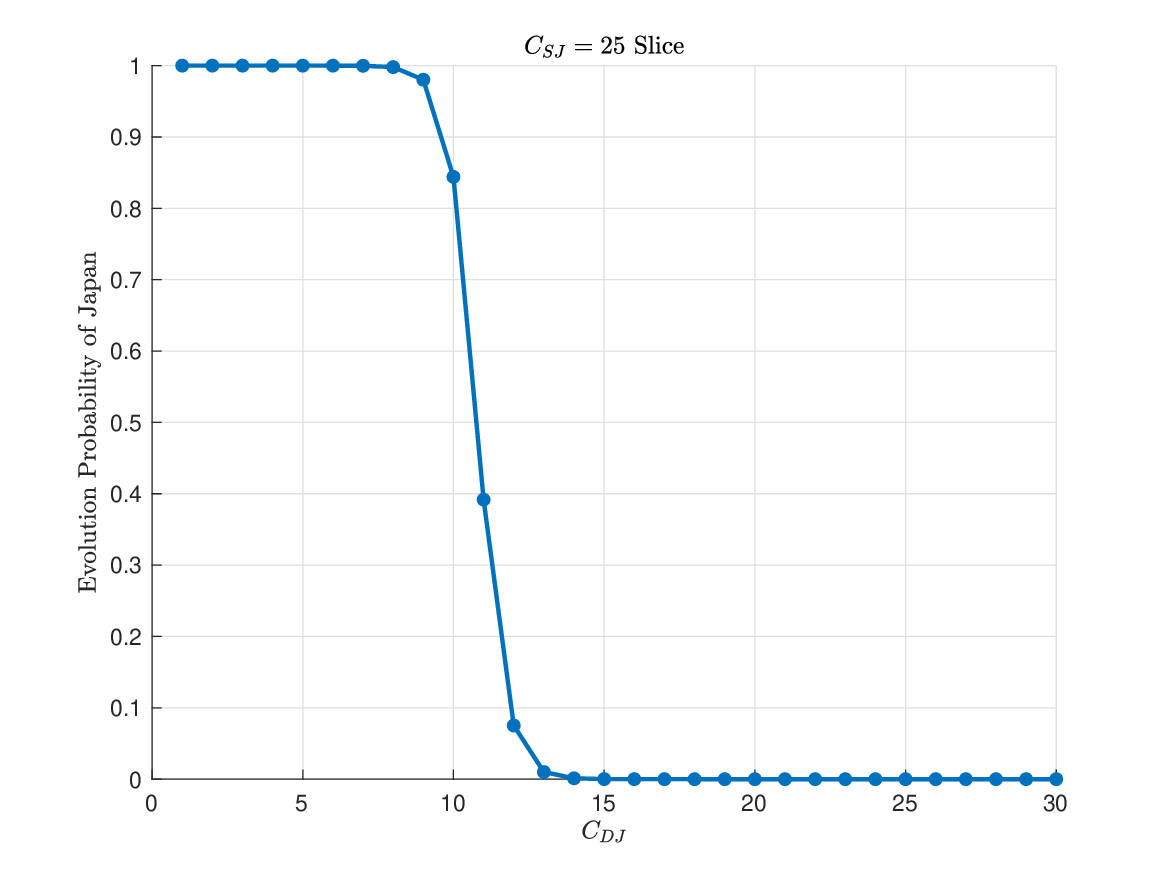}
		\caption{Impact of $C_{DJ}$ on Japan's Evolutionary Decision-making.}
		\label{duo2}
	\end{subfigure}
	\caption{Impact of $C_{SJ}$ and $C_{DJ}$ on Japan's evolutionary decision-making.}
\end{figure}

To analyze the impacts of $C_{SJ}$ and $C_{DJ}$ more precisely, we selected two specific values from Figure~\ref{duo1}, namely $C_{SJ} =25$ and $C_{DJ} =5$, and drew Figures~\ref{duo3} and ~\ref{duo2} respectively. With $C_{DJ} = 5$, as $C_{SJ}$ increases, making the storage costs higher, the Japanese government is more likely to choose the "discharge" strategy. Conversely, with $C_{SJ}=25$, as $C_{DJ}$ increases, to avoid the hefty discharge costs, the Japanese government is more inclined to adopt the "no discharge" strategy.

The above analysis not only emphasizes the central role of cost factors in the Japanese government's decision-making but also suggests that when formulating strategies for handling nuclear wastewater, the government should comprehensively consider the long-term effects on the economy, environment, and society to ensure its decisions are holistic and sustainable.

	\subsubsection{Impact of Japan's nuclear wastewater storage cost $C_{SJ}$ and export tax revenue reduction due to discharge $T_{RJ}$ on Japan's evolutionary decision-making}
	
\begin{figure}[H]
		\centering
		\includegraphics[width=0.6\linewidth]{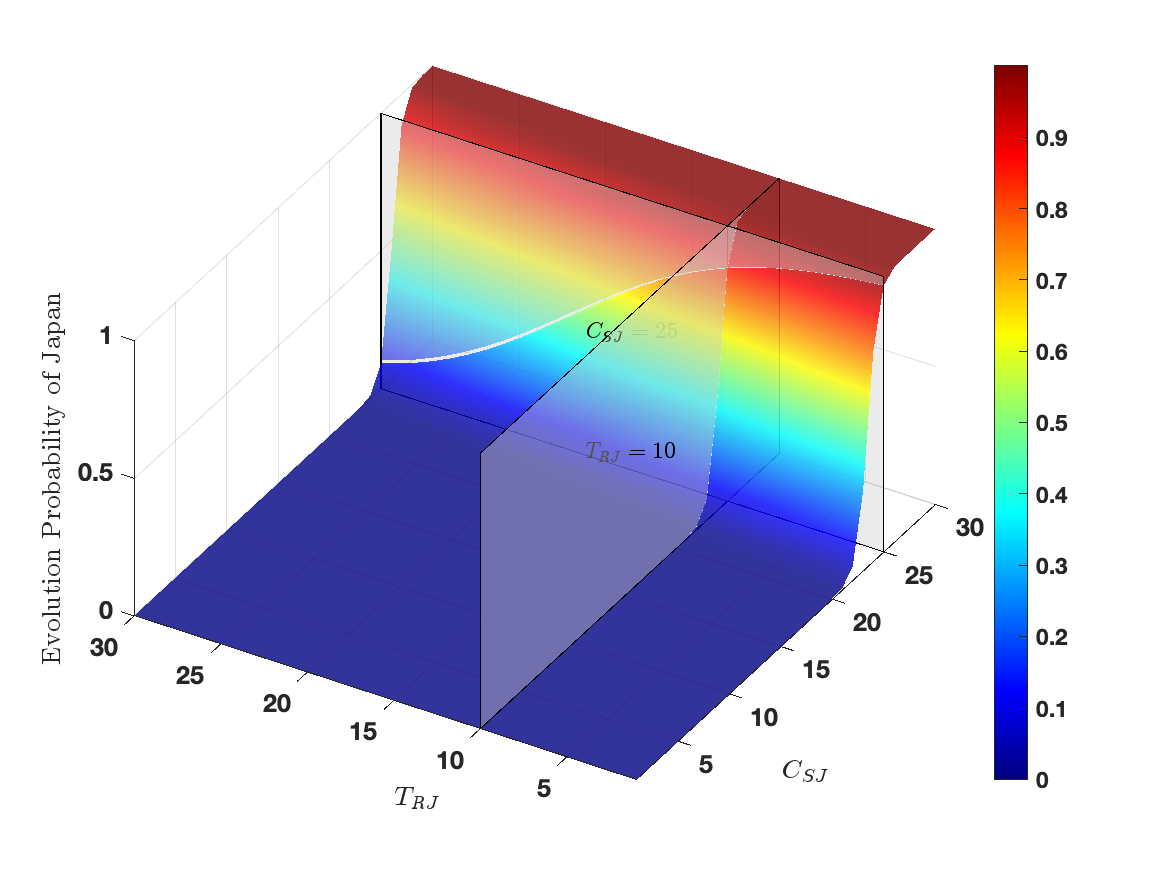}
		\caption{Impact of $C_{SJ}$ and $T_{RJ}$ on Japan's evolutionary decision-making.}
		\label{duo4}
	\end{figure}

	The cost for Japan to store nuclear wastewater, $C_{SJ}$, and the reduction in export tax revenue due to the wastewater discharge, $T_{RJ}$, have a significant impact on the decision-making of the Japanese government. Against this backdrop, to comprehensively understand how $C_{SJ}$ and $T_{RJ}$ influence Japan's decision-making, we base our analysis on the initial parameter settings of Condition 5 from Table~\ref{biao5}. A three-dimensional surface plot was created to provide a visual representation of how these two variables influence Japan's policy stance, as shown in Figure~\ref{duo4}.
	
	Figure~\ref{duo4} clearly reveals the following trends: 
	As $C_{SJ}$ increases, indicating rising costs for storing nuclear wastewater, Japan's policy stance gradually shifts from "no discharge" to "discharge." This is because discharging becomes a more economically viable option as storage costs rise. As $T_{RJ}$ increases, representing the loss in export tax revenue due to the discharge, Japan's stance shifts from "discharge" back to "no discharge." This suggests that when Japan's national revenue declines under a discharge scenario, the government re-evaluates its decision, leaning towards the "no discharge" policy. Notably, under the mentioned parameter settings, when $C_{SJ}$ is less than 20, Japan always opts for the "no discharge" strategy.
	
\begin{figure}[H]
		\centering
		\begin{subfigure}{0.45\textwidth}
			\centering
			\includegraphics[width=1\linewidth]{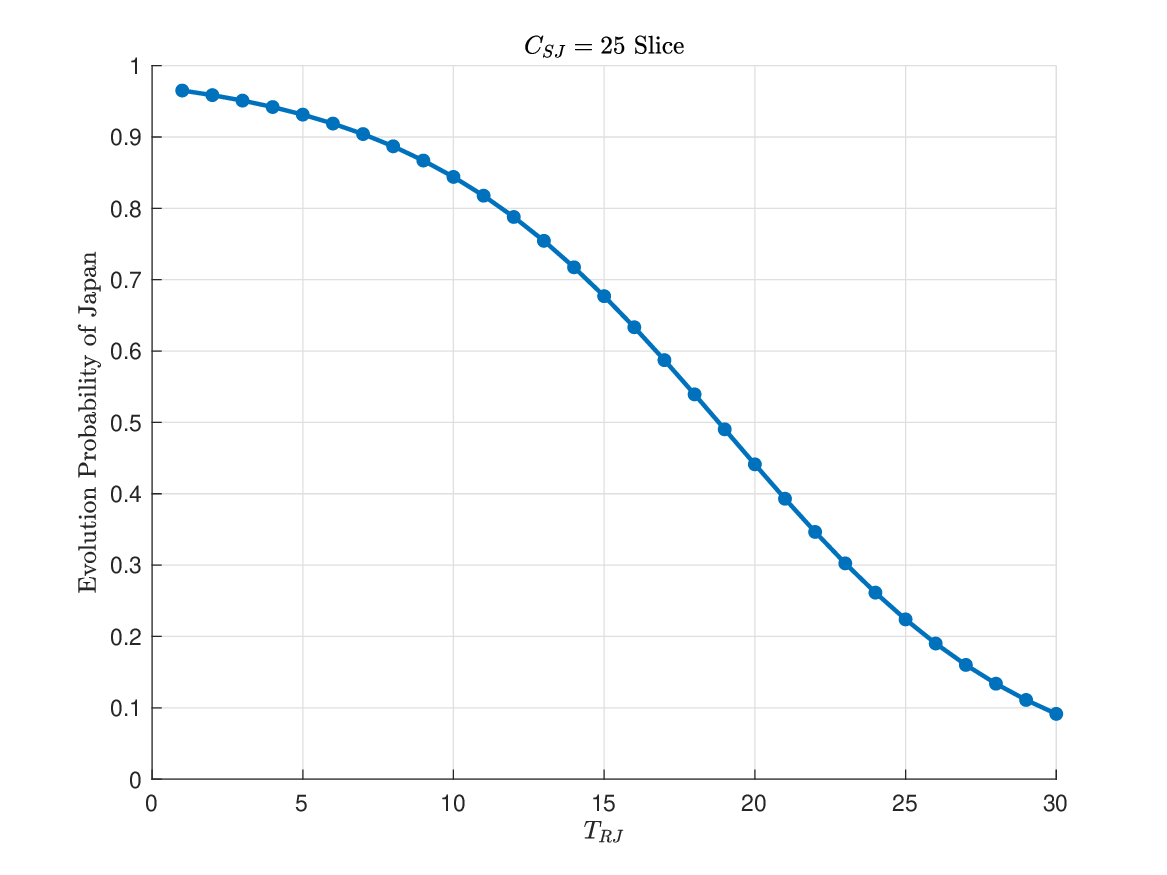}
			\caption{Impact of $C_{SJ}$ on Japan's evolutionary decision-making.}
			\label{duo5}
		\end{subfigure}%
		\begin{subfigure}{0.45\textwidth}
			\centering
			\includegraphics[width=1\linewidth]{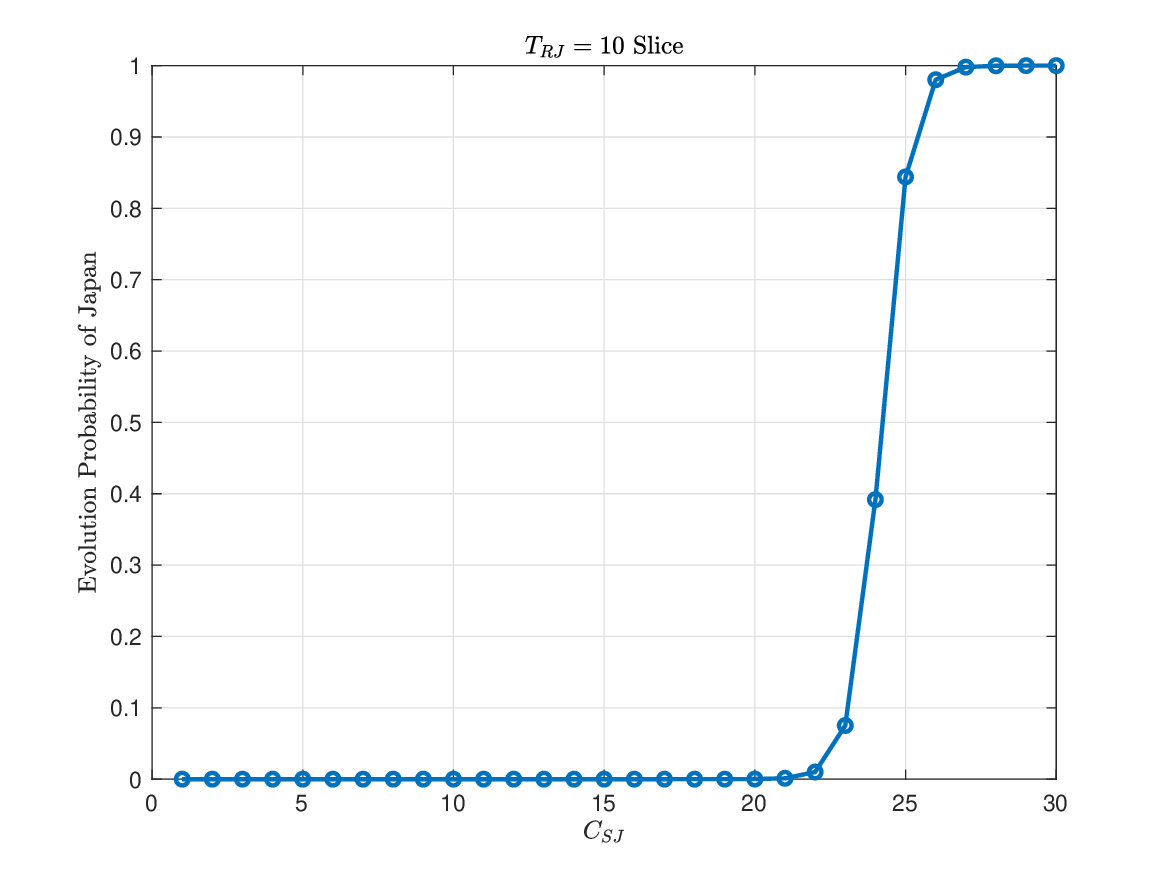}
			\caption{Impact of $C_{HJ}$ on Japan's evolutionary decision-making.}
			\label{duo6}
		\end{subfigure}
		\caption{Impact of $C_{SJ}$ and $C_{HJ}$ on Japan's evolutionary decision-making.}
	\end{figure}

	For a more precise analysis of the impacts of $C_{SJ}$ and $T_{RJ}$, we selected two specific values from Figure~\ref{duo4}, $C_{SJ} =25$ and $T_{RJ} =10$, and drew Figures~\ref{duo5} and ~\ref{duo6} respectively. With $C_{SJ} =25$, as $T_{RJ}$ increases, the loss in export tax revenue due to the discharge policy adversely affects Japan's interests, making the Japanese government more inclined to "continue storage" rather than discharge. On the other hand, with $T_{RJ}=10$, as $C_{SJ}$ rises, the Japanese government is more likely to opt for the "discharge" strategy. This is because as the cost of storing nuclear wastewater escalates, direct discharge becomes a more economically sensible choice.
	
	The above analysis emphasizes the pivotal role of economic factors in the decision-making process. When formulating policies related to the environment, in addition to the direct environmental and health impacts, it's vital to thoroughly consider the potential economic consequences of such decisions.

\subsubsection{Impact of aid from other countries to Japan $C_{HJ}$ and litigation compensation from Japan $C_{LC}$ on the evolutionary strategy of other countries}

\begin{figure}[H]
	\centering
	\includegraphics[width=0.6\linewidth]{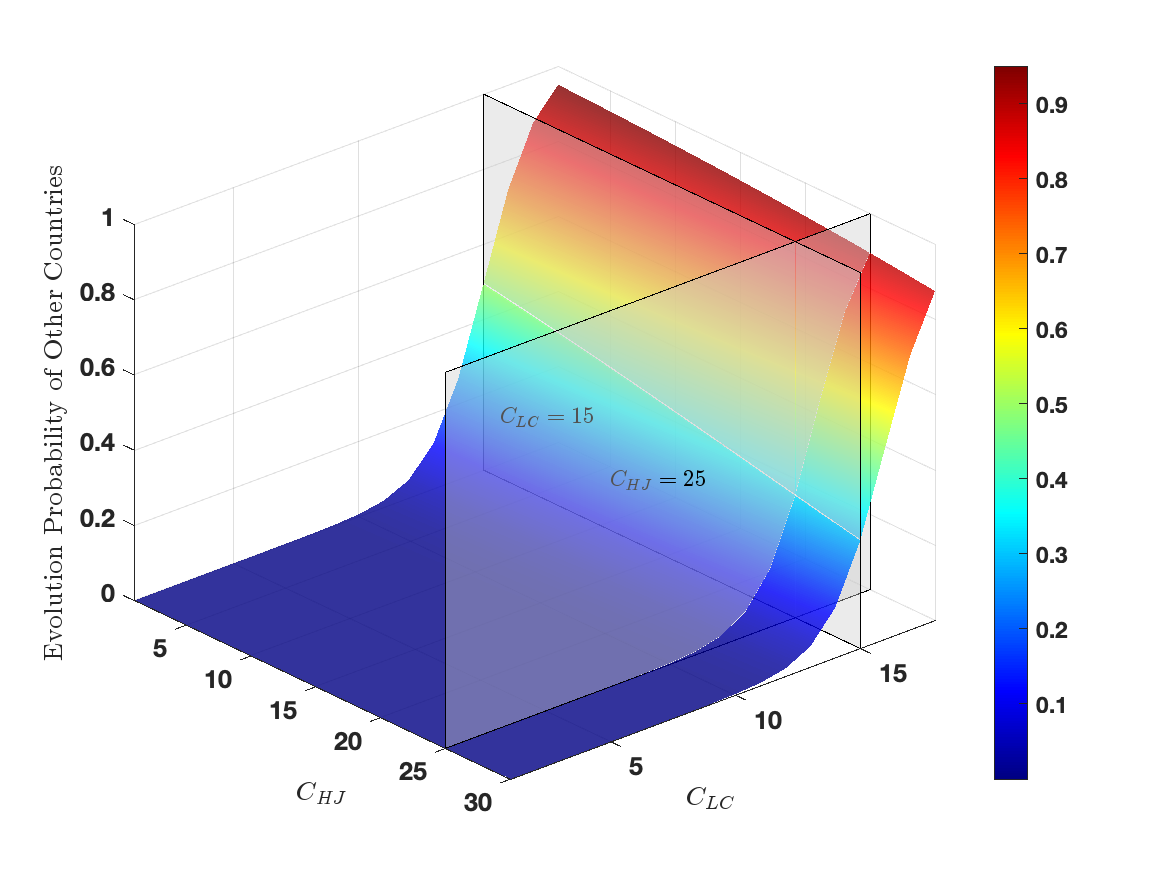}
	\caption{Impact of $C_{HJ}$ and $C_{LC}$ on the evolutionary strategy of other countries.}
	\label{duo7}
\end{figure}

The aid from other countries to Japan, denoted by $C_{HJ}$, and the litigation compensation they receive from Japan, denoted by $C_{LC}$, play pivotal roles in influencing the strategic evolution of these countries. To delve deeper into how these two variables shape the strategic choices of these nations, we refer to the initial parameter settings of Condition 6 in Table~\ref{biao5} and illustrate this relationship with a three-dimensional surface plot as shown in Figure~\ref{duo7}.

Figure~\ref{duo7} reveals how the strategy of other nations towards Japan, specifically their decision to sanction or cooperate, changes with variations in $C_{HJ}$ and $C_{LC}$. As $C_{HJ}$ increases, indicating greater aid to Japan, nations might recognize long-term benefits of strengthened international collaborations, such as technological exchanges and economic partnerships, making them less inclined to sanction Japan. Conversely, as $C_{LC}$ rises, denoting higher litigation compensation, nations might emphasize safeguarding their environmental, economic, and societal interests due to potential damages from Japan's nuclear wastewater discharging actions. This could lead them to favor sanctions as a means to convey their disapproval of Japan's actions.

\begin{figure}[H]
	\centering
	\begin{subfigure}{0.45\textwidth}
		\centering
		\includegraphics[width=1\linewidth]{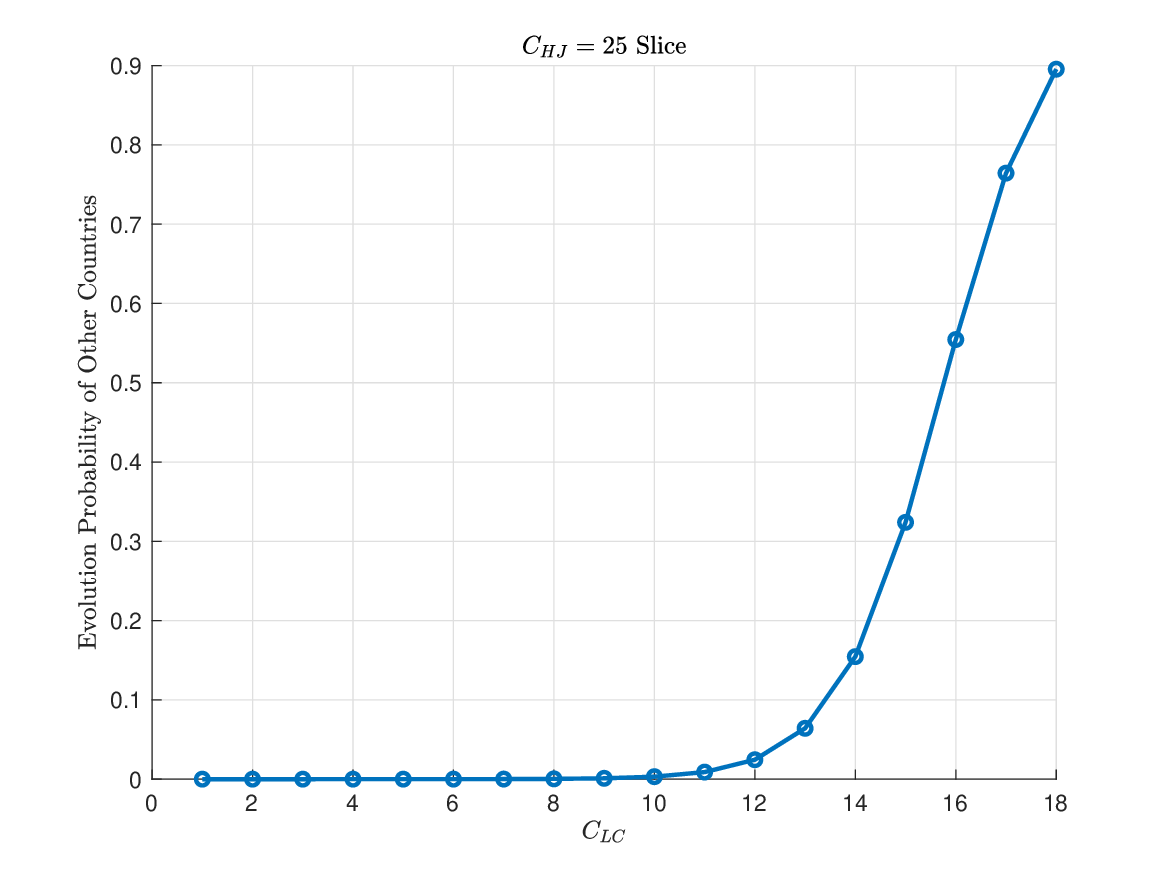}
		\caption{Impact of $C_{LC}$ on the evolutionary strategy of other countries.}
		\label{duo8}
	\end{subfigure}%
	\begin{subfigure}{0.45\textwidth}
		\centering
		\includegraphics[width=1\linewidth]{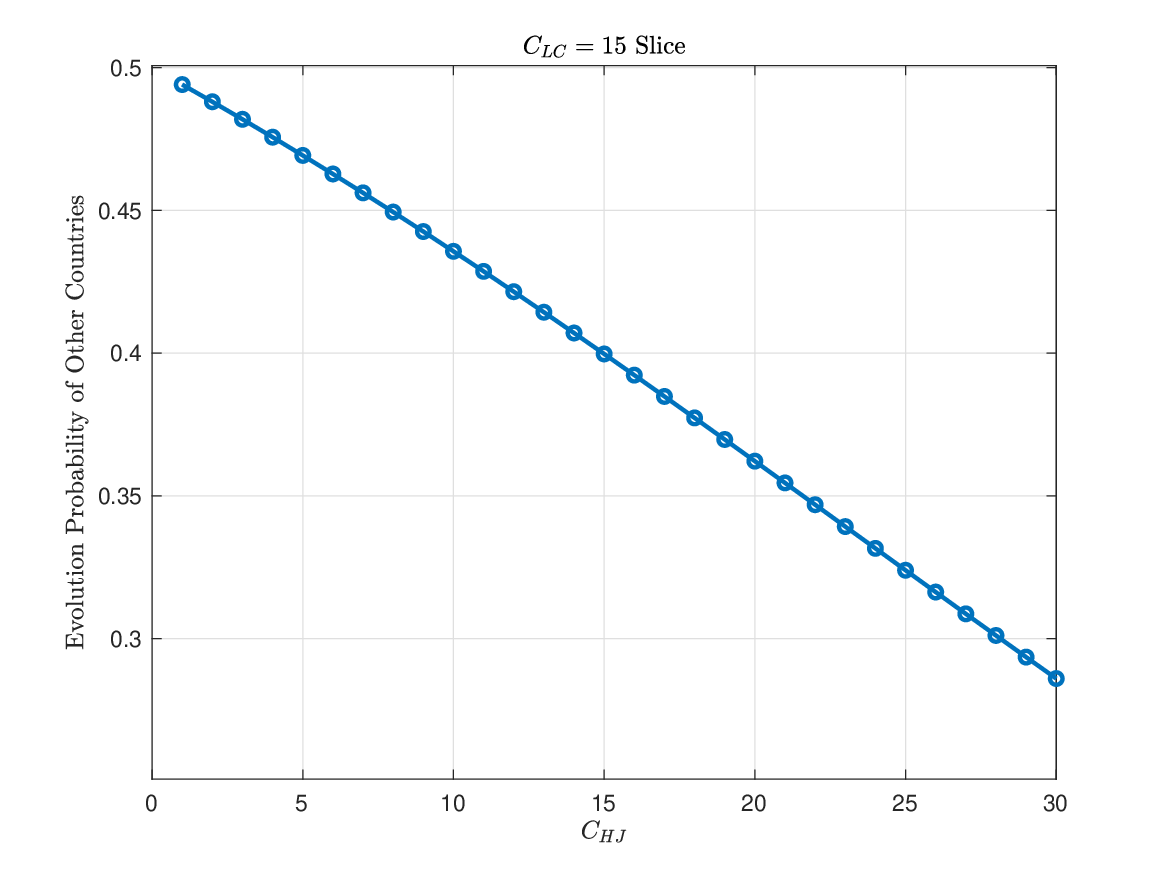}
		\caption{Impact of $C_{HJ}$ on the evolutionary strategy of other countries.}
		\label{duo9}
	\end{subfigure}
	\caption{Influence of $C_{LC}$ and $C_{HJ}$ on the strategic decisions of other nations.}
\end{figure}

To analyze the effects of $C_{HJ}$ and $C_{LC}$ in more detail, we considered specific values $C_{HJ} = 25$ and $C_{LC} = 10$, and plotted them in Figures~\ref{duo8} and~\ref{duo9} respectively. With aid at $C_{HJ} = 25$, nations might be more inclined to sanction Japan as the litigation compensation $C_{LC}$ increases. Conversely, when the litigation compensation is at $C_{LC} = 10$, due to increased aid from other countries to Japan $C_{HJ}$, considering international relations, nations might be more inclined to adopt a non-sanction strategy towards Japan.

The aforementioned analysis highlights the complexities and dynamism inherent in international relations and underscores the intricate balance nations maintain between national interests and global responsibilities. This emphasizes the significance of international collaboration and dialogue in addressing such global challenges.
\section{Conclusions and policy implications}
\subsection{Conclusions}

This study establishes a time-delayed game strategy among Japan, other nations, and the International Atomic Energy Agency (IAEA), presenting a game scenario more closely aligned with real-world situations. By analyzing the characteristic roots of the linearized system, we delve deeply into the stability of the tripartite strategy, and the conditions and possibilities for reaching different asymptotic stable states. We identified the stable points as $\gamma_4(0, 0, 1)$, $\gamma_6(1, 0, 1)$, and $\gamma_8(1, 1, 1)$. To further probe the dynamic properties and evolution of these stable points, we employed replicator dynamic equations and stability conditions to numerically simulate the evolution trajectories of these three points. Finally, we explored the influence of various strategy combinations on the tripartite evolutionary stability strategy, considering Japan's sea discharge cost ($C_{DJ}$), Japan's nuclear wastewater storage cost ($C_{SJ}$), other nations' litigation compensation ($C_{LC}$), aid from other nations to Japan ($C_{HJ}$), and the reduction in export tax revenue due to sea discharge by Japan ($T_{RJ}$).

The key conclusions drawn from this study include:

\begin{itemize}
	\item \textbf{Strategy Stability and Time Delays:}
	
	Introducing time-delay factors for the first time in the context of the evolutionary game of Japan's nuclear wastewater discharge, we observed the evolutionary trajectories under various time delays ($\tau$). Results showed that the impacts of nuclear wastewater on the three parties are not significant in the short term. However, as the time delay gradually increases, there are notable shifts in the strategies of the Japanese government and IAEA. This indicates that time delays have certain effects on strategy stability. Although the final decision results of Japan, other nations, and the IAEA are consistent with those without time delay ($\tau=0$), the presence of time delays means constant adjustments are required during the gaming process. For other countries, based on Condition 1 and Condition 2, the influence of time delays on their evolutionary results is relatively minor. Nevertheless, under Condition 3, with initial values [0.8,0.5,0.5], time delays exerted significant impacts on their evolutionary trajectories for different $\tau$ values. In summary, our study reveals that time delays play a crucial role in the evolution and strategy stability of Japan's nuclear wastewater discharge decisions. Moreover, time delays can also impact the evolutionary trajectories of other countries.
	
	\item \textbf{Urgency of Scientifically Effective Nuclear Wastewater Treatment for Japanese Decision-Making:}
	
	In the analysis of the influence of Japan's nuclear wastewater storage cost $C_{SJ}$ and discharge cost $C_{DJ}$ on its discharge strategy, as the discharge cost $C_{DJ}$ increases, the Japanese government progressively shifts from a "discharge" to a "no-discharge" strategy. Conversely, as the storage cost $C_{SJ}$ rises, the Japanese government transitions from a "no-discharge" to a "discharge" strategy. This cost-driven strategy adjustment underscores the pivotal role of efficient nuclear wastewater treatment technology in Japan's decision-making. Specifically, if Japan can access cutting-edge and cost-effective wastewater treatment technology, associated costs are expected to significantly decrease. This not only provides the Japanese government with broader decision-making options but also aids in ensuring the economic viability of their strategies. Crucially, employing advanced and effective treatment technology can substantially mitigate the environmental risks posed by nuclear wastewater, bolstering public trust and support in government decisions.
	
	\item \textbf{Influence of Reduced Export Tax Revenue on Japan's Nuclear Wastewater Treatment Strategy:}
	
	In the analysis of the influence of Japan's nuclear wastewater storage cost $C_{SJ}$ and reduced export tax revenue due to sea discharge $T_{RJ}$ on its discharge strategy, with a constant storage cost of $C_{SJ}=25$, the Japanese government gradually moves from a "discharge" to a "no-discharge" strategy as $T_{RJ}$ declines. As export tax losses resulting from wastewater discharge increase, the economic pressures on Japan also escalate. In this context, even though the cost of storing nuclear wastewater is relatively high, it becomes a more acceptable expense compared to the loss in export tax revenues. Thus, to safeguard its economic interests and uphold its international image, the Japanese government may lean more towards a "no-discharge" strategy.
	
	\item \textbf{Influence of International Aid and Litigation Compensation on Other Nations' Decisions:}
	
	In the analysis of the influence of litigation compensation to Japan from other nations ($C_{LC}$) and aid to Japan from other nations ($C_{HJ}$) on other nations' sanctioning strategies, as litigation compensation ($C_{LC}$) increases with fixed aid ($C_{HJ}$), other nations are more likely to opt for a strategy to "sanction" Japan, potentially through economic measures such as increasing tariffs or restricting trade. Conversely, with fixed litigation compensation and increasing aid to Japan ($C_{HJ}$), it suggests that other nations see significant benefits in forging closer economic ties with Japan. In this case, other nations are more likely to adopt a "non-sanctioning" strategy towards Japan.
\end{itemize}

	\subsection{Policy implications}
	
	The time-delay game model established between the Japanese government, other countries, and the IAEA in this article, after numerical simulations, clearly shows the strategic outcomes of the three parties. The proposed solution process is also realistic. Based on the research results, we offer the following recommendations:
	\begin{itemize}
		\item 	\textbf{Enhance International Cooperation}
		
The game results show that the strategic evolution of the Japanese government, other countries, and the IAEA is influenced by various factors. This implies that it is difficult for any game participant to achieve the optimal equilibrium state on their own. Therefore, to minimize the impact on their national interests, it is recommended that countries strengthen cooperation with relevant international organizations and openly communicate information and data about marine pollution. This helps relevant countries reduce the loss of benefits brought about by information asymmetry. By sharing data and information, countries can more fully understand the seriousness and potential risks of discharging nuclear wastewater into the sea, and better assess the necessity and impact of taking action.

		\item 	\textbf{Deepen the Study on Time-lag Effects}  
		
This paper first introduces the impact of time-lag effects on the decisions of all parties against the backdrop of nuclear wastewater discharge into the sea. The introduction of time-lag makes the dynamic game model established more realistic and the conclusions more accurate. This type of time-lag evolutionary model provides theoretical references for studying relevant dynamic game issues in other fields. In real life, there may be multiple time-delay factors, such as the effects of ocean currents and other external factors on time delays, and the delay in government decision-making is also an important factor. In practice, the time-lag for each party is not fixed, so future research can delve deeper into the application of distributed time-lag in related problems.

		\item \textbf{Broadly Seek Opinions}  
		
The results of the tripartite game show that the preferences of the parties for strategies may change, emphasizing the necessity for broad communication. Multilateral communication is an important means to promote cooperation and solve the nuclear wastewater issue. Through multilateral communication, countries can directly exchange opinions and interests, find common points of cooperation, and enhance mutual understanding. All parties should listen to each other's opinions and concerns, strive to reach a consensus, and ensure that the interests of all parties are fully considered and balanced in the cooperation process.
		
	\item \textbf{Strengthen Environmental Monitoring and Forecasting}

The time-lag effect of ocean dynamics considered in the model highlights the crucial importance of environmental monitoring in the treatment of nuclear wastewater. To more accurately predict and respond to the potential impact of nuclear wastewater, it is recommended to establish a comprehensive and efficient marine environmental monitoring system to provide a forward-looking assessment of the future transmission path and possible impact areas of nuclear wastewater. In this way, we can not only timely detect and respond to the risks brought by nuclear wastewater but also provide more accurate data support for the model, further enhancing the simulation effect of the game model with distributed time-lag.

	\end{itemize}

	\section{Funding}
	This work is supported by Youth Foundation of China University of Petroleum-Beijing at Karamay  (No.XQZX20230034).
	\section*{Appendix}
\subsection{Proof of Stability Evolution of Different Points}

1. The linearized system corresponding to point $\gamma_{1}(0,0,0)$ is as follows:
\begin{align}
	\frac{d x}{d t} &= \left(C_{S J}-C_{M J}-C_{D J}\right) x(t-\tau) \\
	\frac{d y}{d t} &= -C_{H J} y(t-\tau) \\
	\frac{d z}{d t} &= C_{I I} z(t-\tau)
\end{align}
Let $x=X e^{\lambda t}$, $y=Y e^{\lambda t}$, $z=Z e^{\lambda t}$, and substitute into the linearized system to obtain:
\begin{align}
	\left[\lambda-\left(C_{S J}-C_{M J}-C_{D J}\right) e^{-\lambda \tau}\right] X &= 0 \\
	\left[\lambda+C_{H J} e^{-\lambda \tau}\right] Y &= 0 \\
	\left[\lambda-C_{I I} e^{-\lambda \tau}\right] Z &= 0
\end{align}

The characteristic equation is:
\[
\left[\lambda-\left(C_{S J}-C_{M J}-C_{D J}\right) e^{-\lambda \tau}\right]\left[\lambda+C_{H J} e^{-\lambda \tau}\right]\left[\lambda-C_{I I} e^{-\lambda \tau}\right] = 0
\]

The roots of this transcendental equation satisfy the relationships:

\begin{align}
	\lambda &= \left(C_{S J}-C_{M J}-C_{D J}\right) e^{-\lambda \tau} \\
	\lambda &= C_{I I} e^{-\lambda \tau}  \\
	\lambda &= -C_{H J} e^{-\lambda \tau}
\end{align}

So, this equilibrium point is unstable.

2. The linearized system corresponding to point $\gamma_{2}(1,0,0)$ is as follows:
\begin{align}
	\frac{d x}{d t} &= \left(-C_{S J}+C_{M J}+C_{D J}\right) x(t-\tau) \\
	\frac{d y}{d t} &= \left(-C_{S C}+C_{L C}\right) y(t-\tau) \\
	\frac{d z}{d t} &= C_{I I} z(t-\tau)
\end{align}
Let $x=X e^{\lambda t}$, $y=Y e^{\lambda t}$, $z=Z e^{\lambda t}$, and substitute into the linearized system to obtain:
\begin{align}
	\left[\lambda-\left(-C_{S J}+C_{M J}+C_{D J}\right) e^{-\lambda \tau}\right] X &= 0 \\
	\left[\lambda-\left(-C_{S C}+C_{L C}\right) e^{-\lambda \tau}\right] Y &= 0 \\
	\left[\lambda-C_{I I} e^{-\lambda \tau}\right] Z &= 0
\end{align}

The roots of this transcendental equation satisfy the relationships:

\begin{align}
	\lambda &= \left(-C_{S J}+C_{M J}+C_{D J}\right) e^{-\lambda \tau} \\
	\lambda &= \left(-C_{S C}+C_{L C}\right) e^{-\lambda \tau}   \\
	\lambda &= C_{I I} e^{-\lambda \tau}
\end{align}

So, this equilibrium point is unstable.

3. The linearized system corresponding to point $\gamma_{3}(0,1,0)$ is as follows:
\begin{align}
	\frac{d x}{d t} &= \left(C_{S J}-C_{M J}-C_{D J}-C_{H J}-C_{L C}-I_{J}-T_{R J}\right) x(t-\tau) \\
	\frac{d y}{d t} &= C_{H J} y(t-\tau) \\
	\frac{d z}{d t} &= C_{I I} z(t-\tau)
\end{align}
Let $x=X e^{\lambda t}$, $y=Y e^{\lambda t}$, $z=Z e^{\lambda t}$, and substitute into the linearized system to obtain:
\begin{align}
	\left[\lambda-\left(C_{S J}-C_{M J}-C_{D J}-C_{H J}-C_{L C}-I_{J}-T_{R J}\right) e^{-\lambda \tau}\right] X &= 0 \\
	\left[\lambda-C_{H J} e^{-\lambda \tau}\right] Y &= 0 \\
	\left[\lambda-C_{I I} e^{-\lambda \tau}\right] Z &= 0
\end{align}

The characteristic equation is:
\[
\left[\lambda-\left(C_{S J}-C_{M J}-C_{D J}-C_{H J}-C_{L C}-I_{J}-T_{R J}\right) e^{-\lambda \tau}\right]\left[\lambda-C_{H J} e^{-\lambda \tau}\right]\left[\lambda-C_{I I} e^{-\lambda \tau}\right] = 0
\]

The roots of this transcendental equation satisfy the relationships:

\begin{align}
	\lambda &= \left(C_{S J}-C_{M J}-C_{D J}-C_{H J}-C_{L C}-I_{J}-T_{R J}\right) e^{-\lambda \tau}\\
	\lambda &=C_{H J} e^{-\lambda \tau}  \\
	\lambda &=-C_{I I} e^{-\lambda \tau}
\end{align}

So, this equilibrium point is unstable.

4. The linearized system corresponding to point $\gamma_{4}(0,0,1)$ is as follows:
\begin{align}
	\frac{d x}{d t} &= \left(C_{S J}-C_{M J}-C_{D J}\right) x(t-\tau) \\
	\frac{d y}{d t} &= -C_{H J} y(t-\tau) \\
	\frac{d z}{d t} &= -C_{I I} z(t-\tau)
\end{align}
Let $x=X e^{\lambda t}$, $y=Y e^{\lambda t}$, $z=Z e^{\lambda t}$, and substitute into the linearized system to obtain:
\begin{align}
	\left[\lambda-\left(C_{S J}-C_{M J}-C_{D J}\right) e^{-\lambda \tau}\right] X &= 0 \\
	\left[\lambda+C_{HJ} e^{-\lambda \tau}\right] Y &= 0 \\
	\left[\lambda+C_{I I} e^{-\lambda \tau}\right] Z &= 0
\end{align}

The characteristic equation is:
\[
\left[\lambda-\left(C_{S J}-C_{M J}-C_{D J}\right) e^{-\lambda \tau}\right]\left[\lambda+C_{H J} e^{-\lambda \tau}\right]\left[\lambda+C_{I I} e^{-\lambda \tau}\right] = 0
\]

The roots of this transcendental equation satisfy the relationships:

\begin{align}
	\lambda &= \left(C_{S J}-C_{M J}-C_{D J}\right) e^{-\lambda \tau}\\
	\lambda &=	-C_{H J} e^{-\lambda \tau}  \\
	\lambda &=-C_{I I} e^{-\lambda \tau}
\end{align}

When $C_{S J}<C_{M J}+C_{D J}$ and $C_{I I} < C_{H J}$, this equilibrium point is stable.

5. The linearized system corresponding to point $\gamma_{5}(1,1,0)$ is as follows:
\begin{align}
	\frac{d x}{d t} &= \left(-C_{S J}+C_{M J}+C_{D J}+T_{R J}+C_{H J}+C_{L C}+I_{J}\right) x(t-\tau) \\
	\frac{d y}{d t} &= \left(C_{S C}-C_{L C}\right) y(t-\tau) \\
	\frac{d z}{d t} &= C_{I I} z(t-\tau)
\end{align}
Let $x=X e^{\lambda t}$, $y=Y e^{\lambda t}$, $z=Z e^{\lambda t}$, and substitute into the linearized system to obtain:
\begin{align}
	\left[\lambda-\left(-C_{S J}+C_{M J}+C_{D J}+T_{R J}+C_{H J}+C_{L C}+I_{J}\right) e^{-\lambda \tau}\right] X &= 0 \\
	\left[\lambda-\left(C_{S C}-C_{L C}\right) e^{-\lambda \tau}\right] Y &= 0 \\
	\left[\lambda-C_{I I} e^{-\lambda \tau}\right] Z &= 0
\end{align}

The roots of this transcendental equation satisfy the relationships:

\begin{align}
	\lambda &= \left(-C_{S J}+C_{M J}+C_{D J}+T_{R J}+C_{H J}+C_{L C}+I_{J}\right) e^{-\lambda \tau} \\
	\lambda &= \left(C_{S C}-C_{L C}\right) e^{-\lambda \tau} \\
	\lambda &= C_{I I} e^{-\lambda \tau}
\end{align}

So, this equilibrium point is unstable.

6. The linearized system corresponding to point $\gamma_{6}(1,0,1)$ is as follows:
\begin{align}
	\frac{d x}{d t} &= \left(-C_{S J}+C_{M J}+C_{D J}\right) x(t-\tau) \\
	\frac{d y}{d t} &= \left(C_{L C}-C_{S C}\right) y(t-\tau) \\
	\frac{d z}{d t} &= -C_{I I} z(t-\tau)
\end{align}
Let $x=X e^{\lambda t}$, $y=Y e^{\lambda t}$, $z=Z e^{\lambda t}$, and substitute into the linearized system to obtain:
\begin{align}
	{\left[\lambda-\left(-C_{S J}+C_{M J}+C_{D J}\right) e^{-\lambda \tau}\right] X=0} \\
	{\left[\lambda-(C_{L C}- C_{S C})e^{-\lambda \tau}\right] Y=0} \\
	{\left[\lambda+C_{I I} e^{-\lambda \tau}\right] Z=0}
\end{align}

The roots of this transcendental equation satisfy the relationships:

\begin{align}
	\lambda &=\left(-C_{S J}+C_{M J}+C_{D J}\right) e^{-\lambda \tau} \\
	\lambda &=	\left(C_{L C}-C_{S C}\right) e^{-\lambda\tau} \\
	\lambda &=-C_{I I} e^{-\lambda \tau}
\end{align}

When $C_{S J}>C_{M J}+C_{D J}$ and $C_{L C}<C_{S C}$, this equilibrium point is stable.

7. The linearized system corresponding to point $\gamma_{7}(0,1,1)$ is as follows:
\begin{align}
	\frac{d x}{d t} &= \left(C_{S J}-C_{M J}-C_{D J}-C_{H J}-C_{L C}-I_{J}-T_{R J}\right) x(t-\tau) \\
	\frac{d y}{d t} &= C_{H J} y(t-\tau) \\
	\frac{d z}{d t} &= -C_{I I} z(t-\tau)
\end{align}
Let $x=X e^{\lambda t}$, $y=Y e^{\lambda t}$, $z=Z e^{\lambda t}$, and substitute into the linearized system to obtain:
\begin{align}
	\left[\lambda-\left(C_{S J}-C_{M J}-C_{D J}-C_{H J}-C_{L C}-I_{J}-T_{R J}\right) e^{-\lambda \tau}\right] X &= 0 \\
	\left[\lambda-C_{H J} e^{-\lambda \tau}\right] Y &= 0 \\
	\left[\lambda+C_{I I} e^{-\lambda \tau}\right] Z &= 0
\end{align}

The roots of this transcendental equation satisfy the relationships:

\begin{align}
	\lambda &= \left(C_{S J}-C_{M J}-C_{D J}-C_{H J}-C_{L C}-I_{J}-T_{R J}\right) e^{-\lambda \tau} \\
	\lambda &= C_{H J} e^{-\lambda \tau}  \\
	\lambda &= -C_{I I} e^{-\lambda \tau}
\end{align}

So, this equilibrium point is unstable.

8. The linearized system corresponding to point $\gamma_{8}(1,1,1)$ is as follows:
\begin{align}
	\frac{d x}{d t} &= \left(-C_{S J}+C_{M J}+C_{D J}+T_{R J}+C_{H J}+C_{L C}+I_{J}\right) x(t-\tau) \\
	\frac{d y}{d t} &= \left(-C_{L C}+C_{S C}\right) y(t-\tau) \\
	\frac{d z}{d t} &= -C_{I I} z(t-\tau)
\end{align}
Let $x=X e^{\lambda t}$, $y=Y e^{\lambda t}$, $z=Z e^{\lambda t}$, and substitute into the linearized system to obtain:
\begin{align}
	\left[\lambda-\left(-C_{S J}+C_{M J}+C_{D J}+T_{R J}+C_{H J}+C_{L C}+I_{J}\right) e^{-\lambda \tau}\right] X &= 0 \\
	\left[\lambda-\left(-C_{L C}+C_{S C}\right) e^{-\lambda \tau}\right] Y&= 0 \\
	\left[\lambda+C_{I I} e^{-\lambda \tau}\right] Z &= 0
\end{align}

The roots of this transcendental equation satisfy the relationships:

\begin{align}
	\lambda &= \left(-C_{S J}+C_{M J}+C_{D J}+T_{R J}+C_{H J}+C_{L C}+I_{J}\right) e^{-\lambda \tau} \\
	\lambda &= \left(-C_{L C}+C_{S C}\right) e^{-\lambda \tau} \\
	\lambda &= -C_{I I} e^{-\lambda \tau}
\end{align}

When $C_{L C}>C_{S C}$ and $C_{S J}>C_{M J}+C_{D J}+T_{R J}+C_{H J}+C_{L C}+I_{J}$, this equilibrium point is stable.

\section{Competing interests}
Te authors declare no competing interests.

\end{CJK}
\end{document}